\documentclass[a4paper,11pt,mathcal,amsart]{article}
\usepackage{amsfonts,eucal,amssymb,amsmath,amsthm,amscd} 
\usepackage[all]{xypic}
\newtheorem{theorem}{Theorem}[section]
\newtheorem{corollaire}[theorem]{Corollary}
\newtheorem{remarque}[theorem]{Remark}
\newtheorem{proposition}[theorem]{Proposition}
\newtheorem{lemme}[theorem]{Lemma}
\newtheorem{definition}[theorem]{Definition}

\newtheorem{fait}[theorem]{Fact}
\newenvironment{preuve}{\medskip \noindent {\bf Proof: }}
   {$\diamondsuit$ }

\newcommand{\Ad}{\text{Ad }}

\newcommand{\lieo}{\ensuremath{\mathfrak{o}}}

\newcommand{\liep}{\ensuremath{\mathfrak{p}}}

\newcommand{\lieg}{\ensuremath{\mathfrak{g}}}
\newcommand{\lien}{\ensuremath{\mathfrak{n}}}

\newcommand{\Sn}{\ensuremath{{\bf S}^n}}
\newcommand{\R}{\ensuremath{{\bf R}}}
\newcommand{\RP}{\ensuremath{{\bf RP}}}

\newcommand{\NN}{\ensuremath{{\bf N}}}

\newcommand{\hs}{\ensuremath{\hat s}}
\newcommand{\hd}{\ensuremath{\hat \delta}}

\newcommand{\hm}{\ensuremath{{ \hat M}}}

\newcommand{\hn}{\ensuremath{{ \hat N}}}
\newcommand{\hx}{\ensuremath{{ \hat x}}}
\newcommand{\hy}{\ensuremath{{ \hat y}}}
\newcommand{\hz}{\ensuremath{{ \hat z}}}

\newcommand{\ol}{\ensuremath{{ \omega^L}}}
\newcommand{\om}{\ensuremath{{ \omega^M}}}
\newcommand{\on}{\ensuremath{{ \omega^N}}}

\newcommand{\oo}{\ensuremath{{ \mathfrak o}}}
\newcommand{\lss}{\ensuremath{L \setminus \Lambda}}

\newcommand{\hl}{\ensuremath{{ {\hat \Lambda}}}}

\newcommand{\mcm}{\ensuremath{{ {\mathcal M}}}}
\newcommand{\hmcn}{\ensuremath{{ \hat {\mathcal N}}}}
\newcommand{\mcn}{\ensuremath{{ {\mathcal N}}}}

\newcommand{\dhl}{\ensuremath{{ {\mathcal H}^{n-1}(\Lambda)}}}

\newcommand{\hsi}{\ensuremath{{ \hat \sigma}}}
\newcommand{\si}{\ensuremath{{ \sigma}}}

\begin{document}
\pagenumbering{arabic}
\title{Removable and essential singular sets for higher dimensional conformal maps }
\author{Charles Frances\thanks{Charles.Frances@math.u-psud.fr, partially supported by ANR {\it Aspects Conformes de la G\'eom\'etrie}}}
\date{}
\maketitle

\noindent{{\bf Abstract.} In this article, we prove several results about the extension to the boundary of conformal immersions from an open subset $\Omega$ of a Riemannian manifold $L$,  into another  Riemannian manifold $N$ of the same dimension. In dimension $n \geq 3$, and when  the $(n-1)$-dimensional Hausdorff measure of $\partial \Omega$ is zero, we completely classify the cases when $\partial \Omega$ contains essential singular points, showing that $L$ and $N$ are conformally flat and making the link with the theory of Kleinian groups.}

\section{Introduction}
The aim of this paper is to make progress toward the understanding of singular sets for conformal maps between Riemannian manifolds of dimension at least 3.  The general problem  we are considering can be stated very easily:  assume that  $(L,g)$ and $(N,h)$ are two smooth, connected,  Riemannian manifolds of same dimension $n \geq 2$, and assume that we have a smooth immersion $s : \lss \to N$,   from the complementary of a closed subset $\Lambda  \subset L$, to the manifold $N$, which is conformal, namely $s^*h=e^{\varphi}g$ for some smooth function $\varphi$ on $\lss$.    The set $ \Lambda $ is called a  {\it singular set} for the conformal  immersion $s$, and  a data  $s: L \setminus \Lambda \to N$ as above is  refered to as  {\it a conformal singularity}.  A basic question is to understand under which conditions the singular set $\Lambda$ is removable, namely it is possible to extend $s$ ``across"  $\Lambda$. 

The main contribution of the article is an almost complete understanding of the situation when the dimension $n$ is at least $3$, and the $(n-1)$-dimensional Hausdorff measure of $\Lambda$, denoted ${\mathcal H}^{n-1}(\Lambda)$, is zero.  Under those assumptions, our principal result is Theorem \ref{thm.local}, stated in Section \ref{sec.intro.local} below, which  yields a  local classification of {\it essential}  conformal singularities, namely those for which $s:\lss \to N$ does not extend  to a continuous map from $L$ into the one-point compactification of $N$. Theorem \ref{thm.local} implies that such essential singular sets  can only occur when $L$ and $N$ are conformally flat, and moreover $L$ is a Kleinian manifold.  As a consequence, except in very peculiar situations that are completely classified, singular sets with $\dhl=0$ are removable (maybe adding a point at  infinity to $N$, when $N$ is noncompact), and the extended map is still a conformal immersion (see Theorem \ref{thm.extension}).   Finally, under the extra assumption that $L$ is compact and  the $(n-2)$-dimensional Hausdorff measure of $\Lambda$ is zero, we also classify globally essential conformal singularities in Theorem \ref{thm.global}: in this case $L$ and $N$ are both Kleinian manifolds.  

Since conformal immersions are very peculiar instances in the much larger class of quasiregular mappings, it is natural, before describing our results into more details,  to mention the existing theorems about removable sets and boundary behavior of quasiregular maps.  Quasiregular mappings  (see \cite{iwaniec}, \cite{rickman2},\cite{vaisala1} for comprehensive introductions to the subject) are usually presented  as the ``good" higher dimensional generalization of holomorphic functions of one complex variable.  And indeed, classical theorems of function theory, such as Picard's theorem, or Painlev\'e's theorem on removable sets, find analogous statements in the framework of quasiregular mappings (see for instance \cite{rickman1}, \cite{rickman3},  \cite{vaisala2}).  Most of those results, though, only deal with quasiregular mappings between domains of the extended space $\overline{\R}^n$.  Although more recent works (for instance  \cite{bonk}, \cite{pankka}, \cite{pankka2} and \cite{zorich1}, among others) aimed at some generalizations involving broader classes of  target manifolds $N$, they do not help much for    the problem we are considering, except in very peculiar cases.   
 Moreover, let us stress that  the tools used in the theory of quasiregular mappings involve elaborate analysis, while the very rigid behavior displayed by conformal immersions in higher dimension allow to settle the problem in the conformal framework by purely geometric arguments.  Actually, we hope that the ideas introduced here will be helpful to study removable and essential singular sets for conformal structures which are not Riemannian, the Lorentz signature being of particular interest, and maybe  for other geometric structures of the same kind, such as Cartan geometries.

\subsection{Extension results}
In all the paper, manifolds and maps between them are assumed to be smooth.

We consider as above a  conformal immersion $s: \lss \to N$, where $(L,g)$ and $(N,h)$ are two connected Riemannian manifolds of dimension $n \geq 3$. The conformal structure on $\lss$ is that induced by $(L,g)$.  We will assume that $\dhl=0$, where ${\mathcal H}^{n-1}$ stands for the $(n-1)$-dimensional Hausdorff measure on $(M,g)$ (we refer to \cite[Chap. 4]{mattila} for basic notions on Hausdorff measures).   In particular, $\lss$ is connected and dense in $L$.  In the sequel, those sets satisfying the condition $\dhl=0$ will be refered to as {\it thin singular sets}. 
The points of a (thin) singular set $ \Lambda $ split naturally  into three categories.  

- The {\it removable singular points} are those $x_{\infty} \in  \Lambda $ at which the map $s$ extends continuously.  In other words, there exists a point $y \in N$, so that  for every sequence $(x_k)$ of $\lss$ converging to $x_{\infty}$, the sequence $s(x_k)$ tends to $y$.  

- {\it The poles} are those points $x_{\infty} \in  \Lambda $ such that for every sequence $(x_k)$ of $\lss$ converging to $x_{\infty}$, the sequence $s(x_k)$ leaves every compact subset of $N$. 

- Finally, the points of $ \Lambda$ which are neither removable, nor poles are {\it essential singular points}.  

One thus gets a partition $\Lambda =\Lambda_{rem} \cup\Lambda_{pole}\cup\Lambda_{ess} $ into removable singular points, poles and essential singular points.  We will say that $\Lambda$ is {\it an essential singular set} as soon as $\Lambda_{ess} \not = \emptyset$.

The results of this article will allow  to determine the structure of those three sets for thin singularities, begining with  $\Lambda_{rem}$. 
\begin{theorem}
\label{thm.extension}
Let $(L,g)$ and $(N,h)$ be two connected $n$-dimensional Riemannian manifolds, $n\geq 3$.  Let $\Lambda \subset L$  be a closed subset such that $\dhl=0$, and $s: \lss \to N$  a conformal immersion.
 Then the  set $\Lambda_ {rem}$ is  open   in $\Lambda$ and $s$ extends to a  conformal immersion $s^{\prime}: L \setminus (\Lambda_ {pole}\cup\Lambda_{ess}) \to N$. 

\end{theorem}


In view of the previous theorem, it will be interesting to find criteria ensuring that $\Lambda_{ess}$ is empty.  It turns out that an injectivity assumption on $s$ is enough for that (compare with the result proved in \cite{vaisala2} for  quasiconformal maps).  

\begin{theorem}
\label{thm.extension.embeddings}
Let $(L,g)$ and $(N,h)$ be two connected $n$-dimensional Riemannian manifolds, $n\geq 3$.  Let $\Lambda \subset L$  be a closed subset such that $\dhl=0$, and $s: \lss \to N$  a conformal embedding. Then $\Lambda_{ess}=\emptyset$, and $s$ extends to a conformal embedding $s^{\prime}: L \setminus \Lambda_ {pole} \to N$.   
\begin{enumerate}
\item When $L$ is compact, then $s^{\prime}: L \setminus \Lambda_ {pole} \to N$ is a conformal diffeomorphism.  
\item When both $L$ and $N$ are compact, $\Lambda_{pole}$ is empty, so that $(L,g)$ and $(N,h)$ are conformally diffeomorphic.
\end{enumerate}
\end{theorem}

Assuming that  $L$ is a compact manifold, Theorem \ref{thm.extension.embeddings} classifies,  {\it all possible conformal embeddings} of the Riemannian manifold $(\lss,g)$ into Riemanniann manifolds of the same dimension. It also gives a unicity result for the conformal compactification of $(\lss,g)$: {\it the only compact Riemannian manifold in which $(\lss,g)$ can be embedded as an open subset is $(L,g)$}.

\subsection{Some examples of  essential singular sets}
\label{sec.examples.essential}


Our next step will be to understand, when it is nonempty,  the set $\Lambda_{ess}$ of essential singular points.  Before this, the reader  might like to see examples of conformal immersions admitting essential singularities.  Since there will play a major role in what follows, we describe now  nice  examples  coming from the theory of Kleinian groups.

\subsubsection{Kleinian groups, domain of discontinuity and limit set}
We consider the sphere $\Sn$ with its conformally flat structure.  The conformal group of $\Sn$ is the M\"obius group $PO(1,n+1)$.
One calls  {\it Kleinian group} a discrete subgroup $\Gamma \subset PO(1,n+1)$  which acts freely properly and discontinuously on some nonempty open subset $\Omega \subset \Sn$ (we refer the reader to \cite[Chap. 2]{apanasov}, \cite[Sec. 3.6, 4.6 and 4.7]{kapovich} and  \cite[sec. 5]{matsumoto} for details on the material below).

Given a Kleinian group $\Gamma$, there exists a maximal open set $\Omega({\Gamma}) \subset \Sn$ on which  the action of $\Gamma$ is proper.  This open set $\Omega({\Gamma})$ is called the {\it domain of discontinuity} of $\Gamma$, and its complement in $\Sn$, denoted $\Lambda({\Gamma})$, is called {\it the limit set of $\Gamma$}.  There are several characterizations of the limit set $\Lambda(\Gamma)
$, but two of them will be of particular interest for our purpose.  Let us consider any point $x \in \Omega({\Gamma})$, and denote $\overline{\Gamma.x}$ the closure of the orbit  $\Gamma.x$ into $\Sn$. Then the limit set $\Lambda({\Gamma})$ coincides with $\overline{\Gamma.x} \setminus \Gamma.x$ (see for instance \cite[Lemma 2.2, p 42]{apanasov}). 

Another useful characterization is as follows: the limit set $\Lambda({\Gamma})$ comprises exactly those points $x  \in \Sn$ at which the family $\{  \gamma\}_{\gamma \in \Gamma}$ fails to be equicontinuous (see \cite[Chap. 5]{matsumoto}). The group $\Gamma$ being assumed to be discrete, we observe that its limit set is empty if and only if  $\Gamma$ is finite. 

If $\Gamma \subset PO(1,n+1)$ is a Kleinian group, and $\Omega \subset \Sn$ is a $\Gamma$-invariant open set on which the action of $\Gamma$ is free and properly discontinuous, then the quotient manifold $N:=\Omega/\Gamma$ is naturally endowed with a conformally flat structure, and the the covering map $\pi: \Omega \to N$ is  conformal. Such a quotient $\Omega/\Gamma$ is called a {\it Kleinian manifold}.   When  the action of $\Gamma$ is free on $\Omega(\Gamma)$, the Kleinian manifold $\Omega(\Gamma)/\Gamma$ will be denoted $M(\Gamma)$. It is then the maximal Kleinian manifold that one can build up thanks to the group $\Gamma$. 

\subsubsection{Essential singular sets of Kleinian type}

Let us now consider $\Gamma \subset PO(1,n+1)$ an {\it infinite} Kleinian group, and $\Omega$ an open subset of $\Sn$ on which $\Gamma$ acts freely properly discontinuously.  Let $N:= \Omega\, \slash  \,\Gamma$ be the associated Kleinian manifold. Observe that because we assumed $\Gamma$ infinite, $\Omega$ is a proper open subset of $\Sn$. Denoting by $\Lambda$ the complement of $\Omega$ in $\Sn$, the covering map $\pi$ yields a  conformal singularity $\pi: \Sn \setminus \Lambda \to N$. The set $\Lambda$ turns out to be  an essential singular set for $\pi$.  To see this, we first observe that because $\Gamma$ acts freely properly discontinuously on $\Omega$, we have $\Lambda(\Gamma) \subset \Lambda$. Actually,  $\Lambda(\Gamma) \subset \Lambda_{ess}$.  Indeed, let $x_{\infty} \in \Lambda(\Gamma)$, and let $y$ and $y^{\prime}$ be two distinct points of $N$. Let $z$ and $z^{\prime}$ in $\Omega$ satisfying $\pi(z)=y$ and $\pi(z^{\prime})=y^{\prime}$.   By the characterization of the limit set described above, there exist  two sequences $(\gamma_n)$ and $(\gamma_n^{\prime})$  in $\Gamma$ such that $x_n:=\gamma_n.y$ and $x_n^{\prime}:=\gamma_n^{\prime}.y$ converge to $x_{\infty}$ (actually, we can choose  $\gamma_n=\gamma_n^{\prime}$). Because $\pi(x_n)=y$ while $\pi(x_n^{\prime})=y^{\prime}$, the point $x_{\infty}$ is neither removable, nor a pole, hence is an essential singular point.   On the other hand, let us consider $x_{\infty} \in  \Lambda$
 which is not a pole.  
 It is easily checked that there must be a sequence $(\gamma_n)$ in $\Gamma$ which is not equicontinuous at $x_{\infty}$, so that $x_{\infty} \in \Lambda(\Gamma)$. In particular, $x_{\infty}$ is an essential singular point. The previous discution shows that $\Lambda_{ess}= \Lambda(\Gamma)$ is not empty, and $ \Lambda=\Lambda_{ess} \bigcup \Lambda_{pole}$.  In other words, we have built a conformal singularity $\pi: \Sn \setminus \Lambda \to N$ with an essential singular set  $\Lambda$.
  The conformal singularities constructed in this way will be refered to as {\it conformal singularities of Kleinian type}.

\subsection{Local classification of thin essential singularities}
\label{sec.intro.local}
We are now coming to the main result of our paper, which  asserts basically that locally, all essential conformal singularities having  $(n-1)$-dimensional Hausdorff measure zero are of Kleinian type.  In particular,  the existence of essential singular points imposes strong restrictions on the geometry: the source manifold must be conformally flat, and the target manifold has to be Kleinian.  It is interesting to notice that this geometric restriction does not appear in dimension two, where all Riemannian manifolds are conformally flat.
 
 To state the main theorem, we will need the following definition.

\begin{definition}[minimal essential singular sets]
Let $(L,g)$ and $(N,h)$ be two Riemannian manifolds of dimension $n \geq 2$, and $s:L \setminus \Lambda \to N$ a conformal singularity.  One says that the singular set $\Lambda$  is minimal essential whenever $\Lambda_{rem}=\emptyset$ and $\Lambda_{ess} \not = \emptyset$.
\end{definition}

In view of Theorem \ref{thm.extension},   when studying thin singular sets which are essential,   we can restrict ourselves to minimal essential ones. 
We can finally  state our main result.


\begin{theorem}
\label{thm.local}
Let $(L,g)$ and $(N,h)$ be two connected $n$-dimensional Riemannian manifolds, $n\geq 3$.  Let $\Lambda \subset L$  be a closed subset, such that $\dhl=0$.  Assume that $s: \lss \to N$ is a conformal immersion for which  $\Lambda$ is a minimal  essential singular set. Then:

\begin{enumerate}

\item There exist an infinite  Kleinian group $\Gamma \subset PO(1,n+1)$, a connected open set $\Omega \subset {\bf S}^n$ on which $\Gamma$ acts freely properly discontinuously, and a conformal diffeomorphism $\psi:  N \to \Omega\, \slash  \,\Gamma$.  

\item For each $x_{\infty} \in \Lambda$, there exist an open neighborhood $U \subset L$ containing $x_{\infty}$, and a conformal diffeomorphism $\varphi : U \to V$, where $V$ is an open subset of ${\bf S}^n$, which makes the following diagram commute
$$\xymatrix{
U \setminus \Lambda \ar@{>}[d] ^s \ar@{>}[rr] ^{\varphi}&  & V \setminus \partial \Omega \ar@{>}[d] ^{\pi}\\
N \ar@{>}[rr] ^{\psi}&  & \Omega \,/\, \Gamma
}
$$

In particular, $\varphi(U \cap \Lambda)=V \cap \partial \Omega$ and $\varphi(U \cap \Lambda_{ess})=V \cap \Lambda(\Gamma)$.
\end{enumerate}
\end{theorem}

\subsection{Consequences of the local classification}
Because Theorem \ref{thm.local} classifies locally all thin conformal singularities admitting essential points, the study of a conformal immersion near an essential singular point reduces to understanding what is going on for singularities of Kleinian type.  We can summarize the results in the following corollary.
\begin{corollaire}
\label{coro.picard}
Let $(L,g)$ and $(N,h)$ be two connected $n$-dimensional Riemannian manifolds, $n\geq 3$.  Let $\Lambda \subset L$  be a closed subset, such that $\dhl=0$.  Assume that $s: \lss \to N$ is a conformal immersion. Then:
\begin{enumerate}
\item The set $\Lambda_{ess}$ is closed.  If it is nonempty, it is either  discrete, or perfect. 
\item If $\Lambda_{pole}$ is nonempty, its closure in $\Lambda$ is the set $\Lambda_{pole} \cup \Lambda_{ess}$.
\item Assume that $\Lambda$ is minimal essential. Then for every $x_{\infty} \in \Lambda_{ess}$ and any neighborhood $U$ of $x_{\infty}$ in $L$, $s(U \setminus \Lambda)=N$.

\item If $\Lambda$ is discrete and contains at least one essential singular point, then $\Lambda_{pole}=\emptyset$ and $(N,h)$ is conformally diffeomorphic to a Euclidean manifold, or a generalized Hopf manifold.
\end{enumerate}
\end{corollaire}

We define generalized Hopf manifolds as quotients of $\R^n \setminus \{ 0 \}$ by an infinite discrete subgroup of conformal transformations.  Topologically, those manifolds are finite quotients of ${\bf S}^1 \times {\bf S}^{n-1}$  (see Section \ref{sec.hopf}).

When the singular set $\Lambda$ is reduced to a point, the third and fourth points of the corollary can be compared to Picard's theorem about the behavior of a meromorphic function in the neighborhood of an isolated essential singularity.  Let us also mention that when $s: \lss \to N$ is merely a quasiconformal immersion, and when $\Lambda = \{ p\}$ is an isolated essential singularity, then V.A Zorich proved in \cite{zorich1} and \cite{zorich2} that $s(U \setminus p)=N$ for every neighborhood $U$, and that up to finite quotient, $N$ is homeomorphic to a product ${\bf R}^k \times {\bf T}^{n-k}$ or ${\bf S}^1 \times {\bf S}^{n-1}$.   Its proof does not imply corollary \ref{coro.picard} in the conformal framework, though (see also \cite[Th 2.1 p 81]{rickman2}, \cite{pankka} for other generalizations of Picard's theorem in the quasiregular setting).

\subsection{Global classification of  essential singularities}

Theorem \ref{thm.local}  describes completely the geometry of the target manifold $N$, for a thin  essential conformal singularity $s: \lss \to N$.  The local geometry of $L$ is also determined, but in full generality, we can not expect to determine $L$ globally.  Now, if we assume that $L$ is compact, and under the stronger assumption that the singular set has $(n-2)$-dimensional Hausdorff measure zero, the singularity $s:\lss \to N$ can be described globally.  This is our last main result.

\begin{theorem}
\label{thm.global}
Let $(L,g)$ and $(N,h)$ be two connected $n$-dimensional Riemannian manifolds, $n\geq 3$.  We assume that $L$ is compact. Let $\Lambda \subset L$  be a closed subset, such that ${\mathcal H}^{n-2}(\Lambda)=0$. Assume that $s: \lss \to N$ is  a conformal immersion for which $\Lambda$ is a  minimal essential singular set. Then:
\begin{enumerate}

\item There exists an infinite  Kleinian group $ \Gamma \subset PO(1,n+1)$,  a connected open subset $\Omega \subset {\bf S}^n$ on which $\Gamma$ acts freely properly discontinuously, and a conformal diffeomorphism $\psi: N \to \Omega \,/\, \Gamma$.

\item There exists a subgroup $\Gamma^{\prime} \subset \Gamma$ with $\Lambda({\Gamma^{\prime}} )\subsetneq\Lambda(\Gamma)$, such that $\Gamma^{\prime}$ acts freely properly discontinuously on $\Omega({\Gamma^{\prime}})$, and a conformal diffeomorphism  $\varphi : L \to M(\Gamma^{\prime})$.

\item  Let us call $s^{\prime} : \Omega \,/\, \Gamma^{\prime} \to \Omega \,/\, \Gamma$ the natural covering map, and let us define the closed subsets  $\Lambda^{\prime}$ and $\Lambda_{ess}^{\prime}$ in $M(\Gamma^{\prime})$  as the quotients $(\partial \Omega \setminus \Lambda({\Gamma^{\prime}}))\,/\,\Gamma^{\prime}$ and $(\Lambda(\Gamma) \setminus \Lambda({\Gamma^{\prime}}))\,/\,\Gamma^{\prime}$. Then the conformal diffeomorphism $\varphi$ can be chosen such that $\varphi(\Lambda)=\Lambda^{\prime}$, $\varphi(\Lambda_{ess})=\Lambda_{ess}^{\prime}$, and the following diagram commutes 

$$\xymatrix{
L \setminus \Lambda \ar@{>}[d] ^s \ar@{>}[rr] ^{\varphi}&  & M(\Gamma^{\prime}) \setminus \Lambda^{\prime} \ar@{>}[d] ^{s^{\prime}}\\
N \ar@{>}[rr] ^{\psi}&  & \Omega \,/\, \Gamma
}
$$
\end{enumerate}
\end{theorem}

We will apply this theorem to get a full description of punctual essential singularities on compact manifolds in Theorem \ref{thm.singularite.ponctuelle}.




\subsection{Organization of the paper}
As we already mentioned it, the tools used in this paper are of geometric nature. Especially, the proofs heavily rely on the interpretation of conformal structures (in dimension $\geq 3$) in terms of Cartan geometries.  The necessary background on this topic, as well as the first technical results,  are introduced in Section \ref{section.cartan}.  They allow to begin the study of conformal singularities in Section \ref{section.structure.removable}.  The main point is to understand the behavior of the $2$-jet of a conformal immersion in the neighborhood of the singular set, as explained in Section \ref{sec.holonomy}.  Theorems \ref{thm.extension} and \ref{thm.extension.embeddings} are proved respectively in Sections \ref{sec.removable.set} and \ref{sec.embeddings}.  In Section \ref{sec.conf.flatness}, we show that thin essential singular sets only occur on conformally flat manifolds, an important step toward Theorem \ref{thm.local}.

Section \ref{sec.conf.flat} reviews some basic results about conformally flat structures.  The reader familiar with this material may skip it, except maybe for Section \ref{sec.cauchy.completion} which deals with the less standard notion of Cauchy completion for  conformally flat structures.  This preparatory work allows to complete the proofs of Theorems \ref{thm.local} and \ref{thm.global} in Sections \ref{sec.proof.local} and \ref{sec.proof.global} respectively.  We conclude the paper with Section \ref{sec.isolated}, which provides a full description of punctual essential singularities on compact Riemannian manifolds.

\section{Conformal structures and Cartan connections}
\label{section.cartan}

Let $(L,g)$ be a Riemannian manifold {\it of dimension $n \geq 3$}.  Let $\hat L$ be the bundle of $2$-jets of orthonormal frames on $\hat L$, and $\pi_L: \hat L \to L$ the bundle map.  The bundle $\hat L$ is a $P$-principal bundle over $M$, where $P$ is the conformal group of the Euclidean space ${\bf R}^n$.  The group $P$ is a semi-direct product  $(\R_+^* \times O(n)) \ltimes \R^n$, where the factor $\R_+^* $ corresponds to homothetic transformations of positive ratio, $O(n)$ is the group of linear orthogonal transformations, and  $\R^n$ is identified with the subgroup of translations. Let $\Sn$ be the $n$-dimensional sphere, and $G:=PO(1,n+1)$ the {\it M\"obius group}, namely the group of conformal transformations of the sphere.   The group $P$ is realized as the subgroup of $G$ fixing a point $\nu \in \Sn$.  We denote by $\lieg:=\oo(1,n+1)$ the Lie algebra of the M\"obius group, and by $\liep \subset \oo(1,n+1)$ the Lie algebra of $P$.  

\subsection{Canonical Cartan connection associated to a conformal structure}
\label{section.canonical.cartan}
It is a fundamental fact, known since Elie Cartan, that under the assumption $n \geq 3$, the conformal class $[g]$ defines  on the  bundle $\hat L$ a {\it unique} normal Cartan connection $\ol$ with values in $\oo(1,n+1)$.  The connection $\ol$ is a $1$-form on $\hm$ with values in the Lie algebra $\oo(1,n+1)$, and satisfying the following properties:
\begin{enumerate}
\item{ For every  $\hx \in \hat L$, $\omega_{\hx}^L : T_{\hx}\hat L \to \oo(1,n+1)$ is an isomorphism of vector spaces. }
\item{For every  $X \in \liep$, the vector field $\hat X$ on $\hat L$ defined by $\hat X(\hx):=\frac{d}{dt}_{|t =0}\;\hx.e^{tX}$,  where $Y \mapsto e^Y$ denotes the exponential map on $PO(1,n+1)$, satisfies $\ol(\hat X)=X$. }
\item{For every  $p \in P$, if $R_p$ denotes the  right action by $p$ on $\hm$, then $(R_p)^*\ol = \Ad p^{-1}\ol$.}
\end{enumerate}
The normality condition is put on the curvature of the connection to ensure uniqueness (see {\cite[Chapter IV]{kobayashi}, and more precisely \cite[Theorem 4.2, p\,135]{kobayashi}, as well as  \cite[Chapter 7]{sharpe}).  The triple $(L,\hat L, \ol)$ will be refered to as {\it the normal Cartan bundle} associated to the conformal structure $(L,g)$.  For the conformally flat model  $\Sn=PO(1,n+1)/P$, the normal Cartan bundle is the M\"obius group $G=PO(1,n+1)$, and the Cartan connection is the Maurer Cartan form $\omega^G$.

Let us observe that if $(L,g)$ and $(N,h)$ are two connected $n$-dimensional Riemannian manifolds, $n \geq 3$, and if $s:(L,g) \to (N,h)$ is a conformal immersion, then $s$ lifts to an immersion $\hat s: \hat L \to \hn$, between the bundles of $2$-jets of orthonormal frames.  Moreover, by unicity of the normal Cartan connection, one must have ${\hat s}^*\on=\ol$.  We say that this lift $\hat s$ is {\it a geometric immersion} from $(\hat L,\ol)$ to $(\hn,\on)$.

\subsection{Exponential map}
 \label{section.exponential}
 On the bundle $\hat L$, the Cartan connection $\ol$ yields an exponential map in the following way. The data of  $u$ in ${\mathfrak o}(1,n+1)$ defines naturally a $\omega^L$-constant vector field  $\hat U$ on $\hat L$ by the  relation $\omega^L(\hat U)=u$.  We call $\phi_u^t$ the local flow  generated on $\hat L$ by the field $\hat U$.  At each $\hx \in \hat L$, let ${\mathcal W}_{\hx} \subset {\mathfrak o}(1,n+1)$ be  the set of vectors $u$ such that  $\phi_u^t$ is defined for   $t \in [0,1]$ at $\hx$.  Then one defines the exponential map at  $\hx$ as follows
  $$\exp({\hx},u):= \phi_{u}^1.\hx,\ \ \forall u \in {\mathcal W}_{\hx}$$

 Using the equivariance properties of the Cartan connection listed above, one shows easily  the following important  equivariance property for the exponential map
  \begin{equation}
  \label{equivariance.property}
  \exp({\hx},u).p^{-1}=\exp({\hx.p^{-1}},(\Ad p).u) 
  \end{equation}  
 for every $ u \in {\mathcal W}_{\hx}, \, \, p \in P.$

\subsection{Injectivity radius.}
 The Lie algebra    $\lieo(1,n+1)$ splits as a sum 
 $$\lien^- \oplus {\bf R} \oplus \lieo(n) \oplus \lien^+$$
  where $\liep= {\bf R} \oplus \lieo(n) \oplus \lien^+$ is the Lie algebra of $P$. The algebra corresponding to the factor $\R$ is a Cartan subalgebra. The two abelian $n$-dimensional subalgebras   $\lien^-$ and $\lien^+$ are the root spaces. They are left invariant by the adjoint action of $\R \oplus \lieo(n)$. A detailed description of this material can be found in \cite[Chapter 7]{sharpe}. As we saw, the group $P$ is a semi-direct product $P=(\R_+^* \times O(n)) \ltimes \R^n$. We put on ${\mathfrak o}(1,n+1)$ a scalar product $< \ , \ >$ which is $\Ad O(n)$-invariant, and denote by $||.||$ the norm it induces on ${\mathfrak o}(1,n+1)$.  For every $\lambda>0$, we will denote $B_{\lien^-}(\lambda)$ (resp. $\overline{B}_{\lien^-}(\lambda))$ the open (resp. closed) ball of center $0$ and radius $\lambda$ in $\lien^-$, for the norm $||.||$.  
 The map $u \mapsto \exp({\hx},u)$ is a diffeomorphism from a  sufficiently small neighborhood of $0 \in {\mathfrak o}(1,n+1)$ on its image.  Notice also that because $(\omega_{\hx}^L)^{-1}(\lien^-)$ is transverse to $T_{\hx}(\pi_L^{-1}(x))=(\omega_{\hx}^L)^{-1}(\liep)$, the map $u \mapsto\pi_L \circ  \exp(\hx,u)$ is a diffeomorphism from a sufficiently small neighborhood of $0$ in $\lien^-$, on its image.  We can then define the {\it injectivity radius at $\hx$} as
 $$ \text{inj}_L(\hx):=\inf\{ \lambda>0 \  | \ u \mapsto\pi_L \circ  \exp(\hx,u) \text{ defines an embedding on } B_{\lien^-}(\lambda) \}$$
By the above remarks, $\text{inj}_L(\hx)>0$, and actually $\text{inj}_L(\hx)$ is bounded from below on compact subsets of $\hat L$.
\subsection{Conformal balls, conformal cones}
\label{section.concormal.cones}

We stick to the notations introduced above. Let $S_{\lien^-}$ be the unit sphere of $\lien^-$, with respect to the norm $||.||$. Let ${\cal F}$ be a subset  of $S_{\lien^-}$. In $\lien^-$, we define the cone over ${\cal F}$ of radius $\lambda>0$ as 
$$ {\cal C}({\cal F}, \lambda) = \{ v \in \lien^- \ | \ v=tw \  t \in [0,\lambda], \   w \in {\cal F}   \}$$  
For $x \in L$, $\hx \in \hat L$ in the fiber of $x$,  $0<\lambda<\text{inj}_L(\hx)$, and ${\mathcal F} \subset S_{\lien^-}$, we can define:

- $B_{\hx}(\lambda):=\pi_L \circ \exp(\hx,B_{\lien^-}(\lambda))$, a {\it conformal ball at $x$}.

-  $C_{\hx}({\mathcal F}, \lambda):=\pi_L \circ \exp(\hx,{\cal C}({\cal F}, \lambda)$, a {\it conformal cone of vertex $x$}.

In the model space, namely the standard $n$-sphere $\Sn=PO(1,n+1)/P$, we will simply consider conformal cones with vertex $\nu$, defined by
$$ C({\mathcal F},\lambda):=\pi_G \circ \exp_G( {\cal C}({\cal F}, \lambda))$$
where $\pi_G: PO(1,n+1) \to \Sn$ is the bundle map and $\exp_G$ is the exponential map in $G=PO(1,n+1)$.

Of course, a conformal immersion $s:L \to N$  maps conformal balls/cones of $L$ to conformal balls/cones of $N$. Indeed, it is straigthforward to check the relation 
\begin{equation}
\label{equation-cones}
s(C_{\hx_k}({\cal F}, \lambda))=C_{\hs(\hx_k)}({\cal F}, \lambda)
\end{equation}

%



Our first technical lemma says that it is possible to include ``thick" conformal cones in the complementary of closed sets of $(n-1)$-dimensional Hausdorff measure zero. 
\begin{lemme}
\label{lem.cone.codimension}
Let $(L,g)$ be a Riemannian manifold of dimension $n \geq 3$. Let $\Lambda \subset L$ be a closed subset such that $\dhl=0$. 
For every $\hx \in \hm$, and for every $0<\lambda< \text{inj}_M(
\hx)$, there exists a $G_{\delta}$-dense set ${\cal U}_{\hx} \subset S_{\lien^-}$ such that $C_{\hx}({\cal U}_{\hx},\lambda) \subset \lss$.
\end{lemme}

\begin{preuve}
Let $\hl$ be the inverse image of $\Lambda$ by the bundle map $\pi_L: \hat L \to L$. Let us call $F$ the subset of $\overline{B}_{\lien^-}(\lambda)$ such that $\exp(\hx,F)=\exp(\hx,\overline{B}_{\lien^-}(\lambda)) \cap \hat{\Lambda}$.  By assumption, this set $F$ has $(n-1)$-dimensional Hausdorff measure zero. Let $m_0$ be an integer such that $\frac{1}{m_0} \leq \lambda$.  For every $m \geq m_0$, we call $\pi_{m}: u \mapsto \frac{u}{||u||}$ the radial projection from $A_{m}={\overline B}_{\lien^-}(\lambda) \setminus { B_{\lien^-}}(\frac{1}{m})$ to $S_{\lien^-}$. This is a Lipschitz map, which is moreover closed. Hence, the set $\pi_{m}(F \cap A_m)$ is a closed subset of $S_{\lien^-}$, the  $(n-1)$-dimensional Hausdorff measure of which is zero.  In particular, its complementary ${\cal U}_{m}$ is open and dense in $S_{\lien^-}$.   Thus $ \bigcap_{m\geq m_0}{\cal U}_{m} $ is a $G_{\delta}$-dense subset of $S_{\lien^-}$ that we call ${\cal U}_{\hx}$.  It is now clear by construction that  $C_{\hx}({\cal U}_{\hx},\lambda) \subset \lss$. 
 \end{preuve}

\subsection{Degeneration of conformal cones}
\label{section.degeneracy}
Our aim now is to understand how the ``shape" of a sequence of conformal cones $C_{\hz_k}({\mathcal F},\lambda)$ evolves, as $\hz_k$ leaves every compact subset in $\hat L$.  The answer is partly contained in the lemma below.

\begin{lemme}
\label{to-zero}
Let $(L,g)$ be a Riemannian manifold of dimension $\geq 3$ and $(\hat L,\ol)$ the normal Cartan bundle associated to the conformal structure of $g$. Let $(z_k)$ be a sequence of $L$ converging to $z_{\infty} \in L$.  Let $(\hz_k)$ and  $(\hz_k^{\prime})$ be two lifts of $(z_k)$ in $\hat L$. We assume that $\hz_k$ converges in $\hat L$, while $\hz_k^{\prime}=\hz_k.p_k$ for a sequence $(p_k)$ of $P$ tending to infinity.  Assume that $\inf_{k \in \NN}(\text{inj}_L(\hz_k^{\prime}))>0$. Then for every $0<\lambda<\inf_{k \in \NN}(\text{inj}_L(\hz_k),\text{inj}_L(\hz_k^{\prime}))$, and every ${\cal F} \subset S_{\lien^-}$ such that  $p_k.C({\mathcal F},\lambda) \to \nu$, as $k \to \infty$, for the Hausdorff topology on $\Sn$, we must have $C_{\hz_k^{\prime}}({\mathcal F}, \lambda) \to z_{\infty}$ for the Hausdorff topology on $L$. 
\end{lemme}

\begin{preuve}
This lemma is a particular case of \cite[Lemma 7]{frances1},  (see also \cite[Corollary 3.3]{frances2}), and the reader will find a complete proof there.  The proof involves the notion of developpement of curves, that we don't introduce here.  The upshot is that a conformal cone is a union of conformal geodesics, namely curves of the form $t \mapsto \pi_L \circ \exp(\hx,tu)$, for $u \in \lien^-$. A point $\hx$ in the fiber of $x $ being chosen, one can develop any conformal geodesic passing through $x$ into the sphere $\Sn$, and thus any conformal cone can be developped. For instance, in the situation of Lemma \ref{to-zero}, the developpement of $C_{\hz_k^{\prime}}({\mathcal F},\lambda) $ with respect to $\hz_k$ is $p_k.C({\mathcal F},\lambda)$. Now, the lemma  follows from the fact that conformal geodesics developping on short curves in $\Sn$ are themselves short (\cite[Lemma 3.1]{frances2}), and that  conformal geodesics of $\Sn$ which are Hausdorff-close to $\nu$ must be short (\cite[Proposition 3.2]{frances2}).  \end{preuve}

\section{Extension results}
\label{section.structure.removable}
We consider $(L,g)$ and $(N,h)$ two connected $n$-dimensional Riemannian manifolds, $n \geq 3$.  Let $\Lambda \subsetneq L$ be a closed subset such that $\dhl=0$, and $s: L \setminus \Lambda \to N$  a conformal immersion.  We denote by $(L,\hat L, \omega^L)$ and $(N,\hn,\on)$ the normal Cartan bundles associated to the respective conformal structures, as introduced in section \ref{section.canonical.cartan}.  If $\hat{\Lambda}$ is the inverse image of $\Lambda$ in $\hat L$, then  $(\hat L \setminus \hat{\Lambda}, \ol)$ is the normal Cartan bundle of $(L \setminus \Lambda,g)$. As we saw in \ref{section.canonical.cartan}, we can lift $s$ to a geometric immersion $\hat s: (\hat L \setminus \hat{\Lambda}, \ol) \to (\hn, \on)$. 

\subsection{Holonomy sequences at a boundary point}
\label{sec.holonomy}
Let us consider $x_{\infty} \in  \Lambda $ which is not a pole for $s$.  It means that there exists $(x_k)$ a sequence of $\lss$ which converges to $x_{\infty}$, and such that $s(x_k)$ converges to $y_{\infty} \in N$.  We will actually get more information working in the bundle $\hat L \setminus \hat{\Lambda}$.  Let $\hx_{\infty} \in \hat{\Lambda}$ in the fiber above $x_{\infty}$, and let $(\hx_k)$ be a sequence of $\hat L \setminus \hat{\Lambda}$ projecting on $(x_k)$ and converging  to $\hx_{\infty}$.  The point is that $\hs(\hx_k)$ may not converge in $\hn$, but there always exists a sequence $(p_k)$ such that $\hs(\hx_k).p_k^{-1}$ does converge to a point $\hy_{\infty} \in \hn$ in the fiber of $y_{\infty}$.  

\begin{definition}[holonomy sequence at $x_{\infty}$]
\label{defi.holonomie}
A sequence $(p_k)$ as above will be called {\it a holonomy sequence at $x_{\infty}$} (associated to $(x_k)$).
\end{definition}
The  holonomy sequence $(p_k)$ just encodes the behavior of the $2$-jets of $s$ along the sequence $(x_k)$. Its study will be, as we shall see, a  major tool in understanding the dynamical behavior of $s$ along  $(x_k)$.  In particular, we will see that for thin singular sets $\Lambda$, removable singularities are characterized by bounded holonomy sequences, while essential ones appear together with unbounded  holonomy sequences.  
%


\subsection{Characterization of  removable points by holo\-nomy, and proof of Theorem \ref{thm.extension}}
\label{sec.removable.set}
Our aim now is to characterize the removable singular points in terms of holonomy sequences. These are the contents of the following theorem, which clearly implies Theorem \ref{thm.extension}.

\begin{theorem}
\label{prop.structure.iness}
Let $(L,g)$ and $(N,h)$ be two connected $n$-dimensional Riemannian manifolds, $n \geq 3$.  Let $\Lambda \subset L$ be a closed subset such that $\dhl=0$, and $s:\lss \to N$ a conformal immersion. Let $x_{\infty}$ be a point of $\Lambda_ {ess} \cup\Lambda_{rem}$.  Then the following statements are equivalent:
\begin{enumerate}
\item{The point $x_{\infty}$ is in $\Lambda_ {rem}$.}
\item{There exists $U_{x_{\infty}}$ an open subset of $L$ containing $x_{\infty}$ such that $s$ extends to a conformal immersion $s_{x_{\infty}}: U_{x_{\infty}} \cup(L \setminus \Lambda) \to N$.}
\item{There is a holonomy sequence of $s$ at $x_{\infty}$ which is bounded in $P$.}
\item{All the holonomy sequences of $s$ at $x_{\infty}$ are bounded in $P$.}
\end{enumerate}
\end{theorem}

\begin{preuve}
It is obvious that point $(2)$ implies point $(1)$, and that point $(4)$ implies point $(3)$.  We  just have to show that $(3)$ implies $(2)$, and that $(1)$ implies $(4)$.  


$\bullet$ {\it $(3)$ implies $(2)$.} 

Our  hypothesis is that there is $\hx_{\infty}$ in the fiber  of $x_{\infty}$, a sequence $(\hx_k)$ in $\hat L \setminus \hat{\Lambda}$ converging to $\hx_{\infty}$, and a bounded sequence $(p_k)$ in $P$ such that $\hs(\hx_k).p_k^{-1}$ is converging in $\hn$. Considering subsequences, we may assume that $(p_k)$ has a limit $p_{\infty} \in P$.  Because $\hs(\hx_k.p_k^{-1})=\hs(\hx_k).p_k^{-1}$, we can assume, replacing  $\hx_{\infty}$ by $\hx.p_{\infty}^{-1}$  and $(\hx_k)$ by $(\hx_k.p_k^{-1})$, that $\hy_k:=\hs(\hx_k)$ is converging to  $\hy_{\infty} \in \hn$.   Because $(\hy_k)$ stays in a compact subset of $\hat L$,  we can find  $k_0 \geq 0$, and $0<\lambda<\min(\text{inj}_M(\hx_{k_0}),\text{inj}_N(\hy_{k_0}))$,  such that  $B_{\hx_{k_0}}(\lambda)
$ and $B_{\hy_{k_0}}(\lambda)
$ contain $x_{\infty}$ and $y_{\infty}$ respectively.  

 Lemma \ref{lem.cone.codimension} implies that there exists a $G_{\delta}$-dense subset ${\cal U} \subset S_{\lien^-}$, such that $C_{x_{k_0}}({\cal U},\lambda) \subset \lss$.  
Let us define $s_{x_{\infty}}^{\prime}: B_{\hx_{\infty}}(\lambda) \to N$ by the formula:
$$ s_{x_{\infty}}^{\prime}(\pi_L\circ \exp(\hx_{k_0},u)):=\pi_N \circ \exp(\hy_{k_0},u), \ \ \forall u \in B_{\lien^-}(\lambda)$$
This is a smooth diffeomorphism from $B_{\hx_{\infty}}(\lambda)$ on its image. On the other hand, because $\hs$ is a lift of $s$ and is a geometric immersion, we get for every $u \in {\cal C}({\cal U},\lambda) $
$$ s (\pi_L \circ \exp({\hx_{k_0}},u))=\pi_N \circ \exp(\hs(\hx_{k_0}),u)=\pi_N \circ \exp(\hy_{k_0},u)$$
In other words, $s$ and $s_{x_{\infty}}^{\prime}$ coincide on $C_{x_{k_0}}({\cal U},\lambda)$, which is dense in $B_{\hx_{\infty}}(\lambda) \setminus \Lambda$, hence they coincide on $B_{\hx_{\infty}}(\lambda) \setminus \Lambda $.  But because $\dhl=0$, $B_{\hx_{\infty}}(\lambda) \setminus \Lambda $ is dense in $B_{\hx_{\infty}}(\lambda)$.  As a consequence $s_{x_{\infty}}^{\prime}$ is a  conformal immersion on $B_{\hx_{\infty}}(\lambda)$.  Finally, the map $s_{x_{\infty}} : B_{\hx_{\infty}}(\lambda) \cup (\lss) \to N$ defined by $s_{x_{\infty}}^{\prime}$ on $B_{\hx_{\infty}}(\lambda)$, and $s$ on $\lss$ is well defined, and is a smooth conformal immersion extending $s$.

$\bullet$  {\it $(1)$ implies $(4)$. }

We pick $x_{\infty}$ in $\Lambda_{rem} \cup \Lambda_{ess}$, and we will prove the implication by contradiction. We first need a technical lemma about the dynamics of sequences of $P$ tending to infinity on the set of conformal cones of the sphere.

\begin{lemme}
\label{dynamique-cones}
Let $(p_k)$ be a sequence of $P$ tending to infinity. Then, considering a subsequence of $(p_k)$ if necessary, we are in one of the following cases:
\begin{enumerate}
\item{ For every $\lambda>0$, and every  ball $ {\cal B} \subset S_{\lien^-}$ (for the metric induced by $||.||$) with nonzero radius, there exists ${\cal B}^{\prime} \subset {\cal B}$ a  subball with nonzero radius, such that as $k \to \infty$, $p_k.C({\cal B}^{\prime}, \lambda)\to \nu$ for the Hausdorff topology.}

\item{There exists  a sequence $(l_k)$ of $P$ converging to $l_{\infty}$, such that $l_kp_k$ stays in the factor $\R_+^* $ of $P=(\R_+^* \times O(n)) \ltimes \R^n$, and $(\Ad l_kp_k)(u)=\frac{1}{\lambda_k}u$ for every $u \in \lien^-$, with $\lim_{k \to \infty}\lambda_k=0$.  }
 \end{enumerate}
\end{lemme}
We postpone the proof of  this lemma to Section \ref{section.appendix}, and derive an  interesting consequence for our purpose.   



\begin{lemme}
\label{lem.pas.cas1}
 Assume  that  $x_{\infty} \in \Lambda_{rem} \cup \Lambda_{ess}$ admits a holonomy sequence $(p_k)$ which is unbounded in $P$. Then we are in the second case of Lemma \ref{dynamique-cones}. \end{lemme}

\begin{preuve}   Assume, for a contradiction, that we are in the first case of Lemma \ref{dynamique-cones}.  We get a ball ${\mathcal B} \subset S_{\lien^-}$ with  nonzero radius, and $\lambda>0$, such that $p_k.C({\mathcal B}, \lambda) \to \nu$.  
Let $\hx_{\infty} \in \hat L$ in the fiber of $x_{\infty}$, and $(\hx_k)$  a sequence of $ \hat L\setminus \hat \Lambda$ converging to $\hx_{\infty}$ such that $\hs (\hx_k).p_k^{-1}$ converges to $\hy_{\infty} \in \hn$.  Let $0<\lambda_0< \inf_{k \geq 0}{\text{inj}}_M(\hx_k)$. Lemma \ref{lem.cone.codimension} implies the existence of a $G_{\delta}$-dense subset ${\mathcal U} \subset {\cal B}$ such that for every $k \geq 0$, the cone $C_{\hx_k}({\cal U}, \lambda_0)$ is included in $\lss$. Because $\text{inj}_N(\hs(\hx_k))=\text{inj}_L(\hx_k)$ is bounded from below by a positive number independent of $k$, and because $p_k.C({\cal U}, \lambda)\to \nu$, we can apply lemma \ref{to-zero} for $z_k:=s(x_k)$, $\hz_k^{\prime}:=\hs(\hx_k)$, and $\hz_k:=\hs(\hx_k).p_k^{-1}$. Together with relation (\ref{equation-cones}), this yields  
\begin{equation}
\label{eq.convergence}
s(C_{\hx_k}({\cal U}, \lambda_0)) \to y_{\infty}\,,\, \text{as } k \to \infty
\end{equation} 
 This is actually impossible. Indeed, because $\lambda_0< \inf_{k \geq 0}{\text{inj}}_M(\hx_k)$, we get that for every $k \geq 0$, the map $u \mapsto \pi_L \circ \exp(\hx_k,u)$ is a diffeomorphism from $B_{\lien^-}(\lambda_0)$ on its image. We deduce that any conformal cone $C_{\hx_k}({\cal B}, \lambda_0)$ has nonempty interior, and actually, all the sets $C_{\hx_k}({\cal B}, \lambda_0)$ contain a common open subset $U \subset L$ for  $k \geq k_0$ large enough.  Then, for every $k \geq 0$, $U_k:=U \cap C_{\hx_k}({{\cal U}}, \lambda_0)$ is a $G_{\delta}$-dense subset of $U \setminus \Lambda$, and the same is true for  $U_{\infty}=\bigcap_{k \geq k_0}U_k$.  From  relation (\ref{eq.convergence}), we get $s(U_{\infty})=y_{\infty}$, which contradicts the fact that $s$ is an immersion, hence locally injective  on $U \setminus \Lambda$.  \end{preuve}

We now prove that the existence of an unbounded holonomy sequence associated to $(\hx_k)$ provides some non-equicontinuity phenomena which forbid $x_{\infty}$ to be in $\Lambda_ {rem}$.  The lemma below reflects this fact. It proves actually more than what we just need for the moment, but  this technical statement will also be useful  later.

\begin{lemme}
\label{lem.pas.equicontinu}
Let $x \in \Lambda$, $t_0>0$, and $\gamma: [0,t_0[ \to \lss$ a smooth curve.  We assume that there exists  a sequence $(t_k)$ of $[0,t_0[$ converging to $t_0$ such that $\gamma(t_k)$ converges to $x$, and  $y_k:=s(\gamma(t_k))$ converges to $y \in N$.  We assume that the holonomy sequence associated to $\gamma(t_k)$ is unbounded. Then, there exists $(t_k^{\prime})$ a sequence of $[0,t_0[$ tending to $t_0$, such that $\gamma(t_k^{\prime})$ converges to $x$, and $s(\gamma(t_k^{\prime}))$ converges to $y^{\prime} \in N$, with $y^{\prime} \not =y$. 
\end{lemme}

\begin{preuve}
We choose $(\hx_k)$ a sequence of $\hat L \setminus \hat{\Lambda}$, which project on $\gamma(t_k)$, and which converges to $\hx \in \hat{\Lambda}$ in the fiber of  $x$.  By hypothesis, there exists a sequence   $(p_k)$, which is unbounded in $P$, such that $\hy_k:=\hs(\hx_k).p_k^{-1}$ converges to $\hy \in \hn$ in the fiber of  $y$.   Lemma \ref{lem.pas.cas1} tells us that replacing $(t_k)$ by a subsequence if necessary (which amounts to consider a subsequence of $(\hx_k)$, and the corresponding subsequence of $(p_k)$), and replacing $\hy_k$ by $\hy_k.l_k^{-1}$ for a sequence $(l_k)$ of $P$ tending to $l_{\infty}$, we may assume that the sequence $(p_k)$ is in the factor $\R_+^*$ of $\R_+^* \times O(n) \subset P$.  The Lemma moreover asserts that $(\Ad p_k)(u)=\frac{1}{\lambda_k}u$ for every $u \in \lien^-$, with $\lim_{k \to \infty}\lambda_k=0$.  

We choose $0<r_0<\frac{1}{2} \min_{k \in \NN \cup \{\infty \}}(\text{inj}_M(\hx_k),\text{inj}_N(\hy_k))$, so that  for every $k \in \NN \cup \{\infty \}$, the maps $\varphi_k: u \mapsto \pi_L \circ \exp(\hx_k,u)$ and $\psi_k: u \mapsto \pi_N \circ \exp(\hy_k,u)$ are well defined, and are diffeomorphisms from $B_{\lien^-}(2r_0)$ to open subsets $U_k$ and $V_k$ of $L$ and $N$ respectively.   For every $k \geq 0$, we define  $F_k:=\varphi_k^{-1}(U_k \cap \Lambda)$.

Lemma \ref{lem.cone.codimension}  ensures the existence of a $G_{\delta}$-dense subset ${\cal U} \subset S_{\lien^-}$, such that for every $k \geq 0$, $C_{\gamma(t_k)}({\cal U},2r_0) \subset \lss$.   For $k \geq k_0$ big enough, we will have $2\lambda_k r_0< 2r_0$, and then, Lemma \ref{lem.cone.codimension} amounts to say that ${\cal C}({\cal U},2\lambda_k r_0) \subset B_{\lien^-}(2\lambda_k r_0) \setminus F_k$.  Then, from  relation (\ref{equivariance.property}), we infer that for every $u \in {\cal C}({\cal U},2\lambda_k r_0)$
\begin{equation}
\label{eq.formule}
 \hs(\exp(\hx_k,u)).p_k^{-1}=\exp(\hy_k, \frac{1}{\lambda_k}u)
 \end{equation}
Observing that for each $k$, ${\cal C}({\cal U},2\lambda_k r_0)
$ is dense in $B_{\lien^-}(2\lambda_k r_0) \setminus F_k$, we deduce that  formula (\ref{eq.formule}) holds actually for every $u \in B_{\lien^-}(2\lambda_k r_0) \setminus F_k$.  

Because $\lambda_k \to 0$, the sequence of conformal balls $B_{\hx_k}(2\lambda_k r_0)=\varphi_k(B_{\lien^-}(2\lambda_k r_0))$ tends to $x_{\infty}$ for the Hausdorff topology on $L$.  This means that choosing $k_0 \geq 0$ large enough, we are sure that for $k \geq k_0$, $\gamma([0,t_0[)$ is not included in $B_{\hx_k}(2\lambda_k r_0)$.  In particular, for every $k \geq k_0$, there exists $u_k \in \lien^-$ with $||u_k||=r_0\lambda_k$, and $t_k^{\prime} \in [0,t_0[$, such that $\phi_k(u_k)=\gamma(t_k^{\prime})$. Considering a subsequence, we may assume that $(\frac{u_k}{\lambda_k})$ converges to $v_{\infty}$.  Because the support of $\gamma$ is in $\lss$, we have $u_k \in B_{\lien^-}(2\lambda_k r_0) \setminus F_k$ for every $k \geq k_0$.  Formula (\ref{eq.formule}) then holds, and projecting on $L$ and $N$, we get
$$ s(\varphi_k(u_k))=\psi_k(\frac{u_k}{\lambda_k})$$
Making $k \to \infty$ yields
$$ \lim_{k \to \infty}s(\gamma(t_k^{\prime})) = \psi_{\infty}(v_{\infty})$$
Because $||v_{\infty}||=r_0$ and $\psi_{\infty}$ is a diffeomorphism from $B_{\lien^-}(2r_0)$ on its image, we get that $y_{\infty}^{\prime}=\psi_{\infty}(v_{\infty})$ is different from $y_{\infty}=\psi_{\infty}(0)$.  Finally, because $\gamma(t_k^{\prime})$ tends to $x_{\infty}$, and $\gamma([0,t_0[) \subset \lss$, we see that the only cluster value of $(t_k^{\prime})$ in $[0,t_0]$ is $t_0$.  Hence $t_k^{\prime} \to t_0$, as desired. \end{preuve}

We can now finish the proof  that point $(1)$ implies point $(4)$.   Let $\gamma : [0,1[ \to \lss$ be a smooth curve such that $\gamma(1-\frac{1}{k})=x_k$ for every $k \geq 1$, where $x_k:=\pi_L(\hx_k)$.  Lemma \ref{lem.pas.equicontinu} ensures the existence of a sequence $(t_k^{\prime})$ tending to $1$, such that $\gamma(t_k^{\prime})$ tends to $x_{\infty}$, and $s(\gamma(t_k^{\prime}))$ tends to $y_{\infty}^{\prime} \not = y_{\infty}$.  This forbids $x_{\infty}$ to be in $\Lambda_{rem}$, and we deduce that the existence of an unbounded holonomy sequence implies $x_{\infty} \in \Lambda_{ess}$.    
\end{preuve}

\subsection{Proof of Theorem \ref{thm.extension.embeddings}}
 \label{sec.embeddings}
 We now want to deduce Theorem  \ref{thm.extension.embeddings} from Theorem \ref{prop.structure.iness},  showing that an injectivity assumption on the conformal immersion $s:\lss \to N$ forces $\Lambda_{ess}$ to be empty. Actually, we will see that near an essential singular point, a conformal immersion is highly noninjective, and to formalize this, it is convenient to use the notion of cluster set. Let  $x_{\infty}$  be a point of the singular set $\Lambda$. The cluster set of $x_{\infty}$ is defined as
$$  \text{Clust}(x_{\infty}):=\{  y \in N \  | \  \exists \,  (x_k) \text{ a sequence in } \lss, \ x_k \to x_{\infty}, \ \text{and } s(x_k) \to y    \}$$
The following proposition identifies the cluster set of an essential singular point. 


\begin{proposition}
\label{prop.plat.injectif}
Let $(L,g)$ and $(N,h)$  be two connected $n$-dimensional Riemannian manifolds, $n \geq 3$. Let $\Lambda \subset L$ be a closed subset such that $\dhl=0$, and $s:\lss \to N$ a conformal immersion. Assume that $\Lambda_ {ess}$ is not empty. Then for every $x_{\infty} \in\Lambda_{ess}$, $\text{Clust}(x_{\infty})=N$. In particular, for every neighborhood $U$ of $x_{\infty}$ in $L$, $s(U \setminus \Lambda)$ is a dense open subset of $N$.
\end{proposition}

Proposition \ref{prop.plat.injectif} will be improved later, since we will deduce from Theorem \ref{thm.local} that if $x_{\infty} \in \Lambda_{ess}$, and if $U$ is a neighborhood of $x_{\infty}$ in $L$, we actually have $s(U \setminus \Lambda)=N$ (see Corollary \ref{coro.picard}).  

\begin{preuve} Let $y_{\infty} \in \text{Clust}(x_{\infty})$.  Let us pick $\hx_{\infty}$ in the fiber of $x_{\infty}$, $(\hx_k)$ a sequence of $\hat L \setminus \hat \Lambda$ converging to $\hx_{\infty}$, and $(p_k)$ a sequence of $P$ such that $\hy_k:=\hs(\hx_k).p_k^{-1}$ tends to a point $\hy_{\infty}$ in the fiber above $y_{\infty}$.  By Theorem \ref{prop.structure.iness}, the sequence  $(p_k)$ is unbounded, and Lemma \ref{lem.pas.cas1} ensures that considering  subsequences, we may assume that $(p_k)$ is contained in the factor $\R_+^*$ of $P=(\R_+^* \times O(n)) \ltimes \R^n$.  Moreover, always by Lemma \ref{lem.pas.cas1}, there exists $(\lambda_k)$ a sequence of $\R_+^*$ converging to $0$ such that for every $\mu>0$

\begin{equation}
\label{converge}
 (\Ad p_k).B_{\lien^-}(\mu \lambda_k) = B_{\lien^-}(\mu)
 \end{equation}
If $\mu$ is chosen smaller than $\min_{k \in \NN \cup \{ \infty \}}(\text{inj}_M(\hx_k),\text{inj}_N(\hy_k))$, the maps $u \mapsto \pi_N \circ \exp({\hy_k},u)$  and $u \mapsto \pi_L \circ \exp({\hx_k},u)$ are well defined and diffeomorphisms from $B_{\lien^-}( \mu \lambda_k)$ on their images  for every $k \in \NN \cup \{ \infty \}$. Lemma \ref{lem.cone.codimension} implies the existence of a $G_{\delta}$-dense subset ${\cal U} \subset {S_{\lien^-}}$, such that $C_{\hx_k}({\cal U}, \mu \lambda_k) \subset \lss$ for every $k \geq 0$.  Relations (\ref{equivariance.property}) and (\ref{converge}) then yield
$$s(C_{\hx_k}({\cal U}, \mu \lambda_k)) = C_{\hy_k}({\cal U}, \mu)$$
In particular, $s(C_{\hx_k}({\cal U}, \mu \lambda_k)) \to C_{\hy_{\infty}}({\cal U},\mu)$ as $k \to \infty$. We infer that $ C_{\hy_{\infty}}({\cal U},\mu) \subset \text{Clust}(x_{\infty})$, and finally $B_{\hy_{\infty}}(\mu) \subset \text{Clust}(x_{\infty})$ because $\text{Clust}(x_{\infty})$ is a closed set.  Since $B_{\hy_{\infty}}(\mu)$ is a neighborhood of $y_{\infty}$, we just showed that $\text{Clust}(x_{\infty})$ is an open set. We assumed that $N$ is connected, so that we get $\text{Clust}(x_{\infty})=N$. In particular, for every neighborhood $U$ of $x_{\infty}$ in $L$, we must have $\overline{s(U \setminus \Lambda)}=N$, hence $s(U \setminus \Lambda)$ is a dense open subset of $N$.  \end{preuve}

We can now prove Theorem \ref{thm.extension.embeddings}.  Observe first that Proposition \ref{prop.plat.injectif} above ensures that if $s: \lss \to N$ admits essential singular points, then $s$ can not be injective. We infer that $\Lambda=\Lambda_{rem} \cup \Lambda_{pole}$. By Theorem \ref{thm.extension}, we know that $L \setminus \Lambda_{pole}$ is an open subset of $L$, and that $s$ extends to a conformal immersion $s^{\prime}: L \setminus \Lambda_{pole} \to N$.  Actually $s^{\prime}$ is injective, hence an embedding.  Indeed, if $s^{\prime}$ is not injective, we can find two disjoint open sets $U$ and $V$ in $L \setminus \Lambda_{pole}$ such that $s^{\prime}$ maps $U$ and $V$ diffeomorphically on the same open set $W$.  Because $s^{\prime}(U \cap (L \setminus \Lambda))$ and  $s^{\prime}(U \cap (L \setminus \Lambda))$ are two dense open subsets of $W$, they intersect, contradicting the injectivity of $s$ on $L \setminus \Lambda$.  

Assuming that $L$ is compact, the definition of poles implies that the immersion $s^{\prime}: L \setminus \Lambda_{pole} \to N$ is a proper map.  By connectedness of $N$, it has to be onto. Finally $s^{\prime}$ is a conformal diffeomorphism between $(L \setminus \Lambda_{pole},g)$ and $(N,h)$.  

If moreover $N$ is also assumed to be compact, then $\Lambda_ {pole}$ is empty, and we get that   $(L,g)$ and $(N,h)$ are conformally diffeomorphic.

\subsection{Essential singular points imply conformal flatness}
\label{sec.conf.flatness}
We are now going to make an important step toward Theorem \ref{thm.local}, proving that the existence of thin essential singular sets is only possible on conformally flat manifolds. Thus, generically,  by Theorem \ref{thm.extension}, if a thin singular set contains no poles (for instance if $N$ is compact), it is always possible to extend a conformal immersion across  it.
In the following, by {\it conformal curvature} on a Riemannian manifold, we will mean the Weyl curvature tensor when the dimension is $\geq 4$, and the Cotton tensor when the dimension is $3$ (see \cite[p\,131]{akivis}). 
\begin{proposition}
\label{prop.essential.flat}
Let $(L,g)$ and $(N,h)$  be two connected $n$-dimensional Riemannian manifolds, $n \geq 3$. Let $\Lambda \subset L$ be a closed subset such that $\dhl=0$, and $s:\lss \to N$ a conformal immersion. Assume that $\Lambda_ {ess}$ is not empty. Then for every $x_{\infty} \in\Lambda_{ess}$, and every $y_{\infty}$ in $\text{Clust}(x_{\infty})$, the conformal curvature vanishes at $y_{\infty}$. In particular, the manifolds $(L,g)$ and $(N,h)$ are both conformally flat.
\end{proposition}

\begin{preuve} We pick $y_{\infty} \in \text{Clust}(x_{\infty})$, and we consider $\hx_{\infty},\hx_k,\hy_k,\hy_{\infty},p_k, \mu$ and ${\mathcal U}$ as at the begining of the proof of Proposition \ref{prop.plat.injectif}.  On $\hat L$, there is, associated to the normal Cartan connection $\omega^L$, a {\it curvature function}  $\kappa$  (we don't give details here, and refer the reader to \cite{sharpe} chapters 5.3 and 7).  This is a map $\kappa: \hat L \to \text{Hom}(\Lambda^2({\mathfrak o}(1,n+1)\,/\,\liep),\liep)$, satisfying the equivariance relation:
\begin{equation}
\label{equivariance.courbure}
\kappa_{\hx}(v,w)= (\Ad p^{-1}).\kappa_{\hx.p^{-1}}((\Ad p).v,(\Ad p).w).
\end{equation}
 The vanishing of the Cartan curvature $\kappa$ at $\hx$ implies the vanishing of $\kappa$ on the fiber of $\hx$. It thus makes sense to say that $\kappa$ vanishes at a point $x \in L$, and this is equivalent to the vanishing of the conformal curvature at $x$ (see \cite[Chap. 7]{sharpe}).  Hence, to get the lemma, it is enough to show that $\kappa$ vanishes at $y_{\infty}$.  
 
For convenience, we will see $\kappa$ as a map from $\hat L$ to  $\text{Hom}(\Lambda^2(\lien^-),\liep)$.  Then, relation (\ref{equivariance.courbure})  still holds, provided $p \in \R_+^* \times O(n) \subset P$. Now, since $\hs$ is a geometric immersion, we have for every $v,w \in \lien^-$, and every $k \in \NN$
$$ \kappa_{\hx_k}(v,w)=\kappa_{\hs(\hx_k)}(v,w)$$
By  relation (\ref{equivariance.courbure}),  we also get
$$ \kappa_{\hs(\hx_k)}(v,w)=(\Ad p_k^{-1}).\kappa_{\hy_k}((\Ad p_k).v,(\Ad p_k).w)$$
Recall that  $\Ad p_k^{-1}$ (resp. $\Ad p_k$) acts trivially on ${\bf R} \oplus {\mathfrak o}(n)$, and by multiplication by $\frac{1}{\lambda_k}$ on $\lien^+$ (resp. $\lien^-$).  Writting $\kappa_{\hy_k}^{(1)}(v,w)$ and $\kappa_{\hy_k}^{(2)}(v,w)$ for the components of $\kappa_{\hy_k}(v,w)$ on ${\bf R} \oplus {\mathfrak o}(n)$ and $\lien^+$ respectively,  the last two equalities yield
$$  \kappa_{\hx_k}(v,w) = \frac{1}{\lambda_k^2}  \kappa_{\hy_k}^{(1)}(v,w) + \frac{1}{\lambda_k^3}  \kappa_{\hy_k}^{(2)}(v,w) $$
Since $\lambda_k \to 0$, making $k \to \infty$ gives $ \kappa_{\hy_{\infty}}(v,w) =0$, and finally $\kappa_{\hy_{\infty}}=0$.  
The conformal curvature vanishes on $\text{Clust}(x_{\infty})$, and by Proposition \ref{prop.plat.injectif}, $\text{Clust}(x_{\infty})=N$, so that $(N,h)$ is conformally flat.    The manifold $(\lss,g)$ is mapped into $(N,h)$ by a conformal immersion, hence $(\lss,g)$ is itself conformally flat.  Finally, because $\dhl=0$, $\lss$ is dense in $L$, and we get that $(L,g)$ is also conformally flat.  \end{preuve}

\section{Background on conformally flat manifolds}
\label{sec.conf.flat}
By Proposition \ref{prop.essential.flat}, conformal singularities $s: L \setminus \Lambda \to N$ such that  $\dhl=0$ and $\Lambda_ {ess} \not =\emptyset$  only occur when $L$ and $N$ are conformally flat.  To prove Theorem \ref{thm.local}, and in particular to get that $N$ has to be a Kleinian manifold, we will need basic notions about conformally flat manifolds.  Good general references on the subject are \cite{goldman}, \cite[Sec. 3]{matsumoto}  and \cite[Chap. 3, p 139]{thurston}.  All manifolds in the sequel are still assumed to have dimension $\geq 3$.

\subsection{Holonomy coverings}

Among conformally flat manifolds, a nice subset comprises those who admit conformal immersions into  the sphere.  Such immersions are called {\it developping maps}.  When it exists, a developping map is essentially unique.

\begin{fait}
\label{fact.post.compone}
If $(M,g)$ is a connected conformally flat manifold of dimension $n \geq 3$, and if $\delta_1$, $\delta_2$ are two conformal immersions from $M$ to ${\bf S}^n$, then there exists an element $g$ of the M\"obius group such that $\delta_2=g \circ \delta_1$.
\end{fait}
The key point to get the fact above is Liouville's theorem (see for instance \cite[p  310]{spivak}):  {\it a conformal immersion between two connected open subsets $U$ and $V$ of ${\bf S}^n$, $n \geq 3$, is the restriction of a M\"obius transformation.}

  One thus get a M\"obius transformation $g$ such that the set where $\delta_2=g \circ \delta_1$ is nonempty and has empty boundary.

Fact \ref{fact.post.compone} easily  implies that if $\delta: M \to {\bf S}^n$ is a developping map, there  exists a group homomorphism 
$$\rho : \text{Conf}( M, [ g]) \to PO(1,n+1),$$ called {\it the holonomy morphism} associated to $\delta$, such that for every $\varphi \in \text{Conf}( M, [ g])$
\begin{equation}
\label{eq.developping.map}
\delta \circ \varphi = \rho(\varphi) \circ \delta
\end{equation}

Let us now consider a conformally flat structure $(M,[g])$. It is a classical result, which already appears in \cite{kuiper} (see also \cite[Sec. 3]{matsumoto}) that the universal covering $(\tilde M, [\tilde g])$, endowed with the lift $[\tilde g]$ of the conformal structure $[g]$, admits a developping map $\tilde \delta: \tilde M \to {\bf S}^n$.  Let us identify $\pi_1(M)$ with a discrete subgroup $\Gamma \subset \text{Conf}( \tilde M, [ \tilde g])$, and call $\Gamma_{\tilde \rho}$ the Kernel of the holonomy morphism $\tilde \rho : \Gamma \to PO(1,n+1)$. The developping map $\tilde \delta$ induces a conformal immersion $\delta$ from the quotient manifold ${\mathcal M}:=\tilde M\,/\,\Gamma_{\tilde \rho}$ to ${\bf S}^n$.  This manifold ${\mathcal M}$ is called the {\it holonomy covering of } $M$. It is in some sense the ``smallest" conformal covering of $M$ admitting a conformal immersion to the sphere.  This is the meaning of the following lemma.

\begin{lemme}
\label{lem.holonomy.minimal}
Let $M$ be a  connected $n$-dimensional conformally flat Riemannian manifold, $n\geq 3$, and ${\mcm}$   the holonomy covering of $M$. Assume that ${\mathcal M}^{\prime}$ is another connected $n$-dimensional conformally flat Riemannian manifold such that:
\begin{enumerate}
\item There exists a conformal immersion $\delta^{\prime}: {\mathcal M}^{\prime} \to {\bf S}^n$.
\item There exists a conformal covering map $\pi: {\mathcal M}^{\prime} \to M$.
\end{enumerate}
Then there exists a conformal covering map from ${\mcm}^{\prime}$ onto $\mcm$. 
\end{lemme} 

\begin{preuve}
Let us call $\tilde M$ the conformal universal covering of $M$, and identify $\pi_1(M)$ with a discrete group $\Gamma$ of conformal transformations of $\tilde M$, so that $M$ is conformally diffeomorphic to $\tilde M \,/\, \Gamma$.   Because $\mcm^{\prime}$ is a covering of $M$, there exists $\Gamma^{\prime} $ a subgroup of $ \Gamma$ such that $\mcm^{\prime}$ is conformally equivalent to $\tilde M\,/\, \Gamma^{\prime}$.  The immersion $\delta^{\prime}$ lifts to a conformal immersion $\tilde{\delta}^{\prime}: \tilde M \to {\bf S}^n$.  Let  $\tilde \delta: \tilde M \to {\bf S}^n$ be a developping map, and $\tilde \rho: \Gamma \to PO(1,n+1)$ the associated holonomy morphism. By Fact \ref{fact.post.compone}, there exists $g \in PO(1,n+1)$ such that $\tilde{\delta}^{\prime}=g \circ \tilde \delta$.  Now, for every $
\gamma \in \Gamma^{\prime}$, one has $\tilde{\delta}^{\prime} \circ \gamma=\tilde{\delta}^{\prime}$, so that $g \circ \tilde \delta \circ \gamma=g \circ \tilde \delta$.  Finally, we get  $\Gamma^{\prime} \subset \Gamma_{\tilde \rho}=\text{Ker }\tilde \rho$.   Hence, there is a conformal covering map from $\mcm^{\prime}=\tilde M\,/\,\Gamma^{\prime}$ onto $\mcm=\tilde M\,/\,\Gamma_{\tilde \rho}$. 
\end{preuve}

\begin{lemme}
\label{lem.lift}
Let $M$ and $N$ be two connected $n$-dimensional conformally flat manifolds, $n \geq 3$.  Let $\mcn$ be the holonomy covering of $N$. Assume there exists  a conformal immersion $\delta: M \to {\bf S}^n$.  Then any conformal  immersion $s:M\to N$ can be  lifted to a conformal immersion $\sigma: M \to \mcn$.\end{lemme}

\begin{preuve}
Let $\tilde M$ and $\tilde N$ be the conformal universal coverings of $M$ and $N$ respectively, and $\pi_M: \tilde M \to M$, $\pi_N: \tilde N \to N$ the associated  covering maps.  We denote by $\Gamma_{M}$ and $\Gamma_{N}$ the fundamental groups of $M$ and $N$, seen as discrete subgroups of conformal transformations of $\tilde M$  and $
\tilde N$.  The conformal immersion $\delta$ lifts to a developping map $\delta_M:\tilde M \to {\bf S}^n$, satisfying $\delta_M \circ \gamma=\delta_M$ for every $\gamma \in \Gamma_M$. We also introduce  $\delta_N$ a  developping map on $\tilde N$, and denote by  $\rho_N: \Gamma_N \to PO(1,n+1)$ the associated holonomy morphism.  The conformal immersion $s$ lifts to a conformal immersion $\tilde s: \tilde M \to \tilde N$, and there is a morphism $\rho: \Gamma_M \to \Gamma_N$ such that for every $\gamma \in \Gamma_M$, $\tilde s \circ \gamma= \rho(\gamma) \circ \tilde s$.  Thanks to Fact \ref{fact.post.compone},   there exists an element $g \in PO(1,n+1)$ such that $ \delta_N \circ \tilde s=g \circ \delta_M$.  For every $x \in \tilde M$ and every $\gamma \in \Gamma_M$, we have on the one hand
$$ \delta_N(\rho(\gamma).\tilde s(x))=\rho_N(\rho(\gamma)).\delta_N(\tilde s(x))$$  
and on the other hand
$$\delta_N(\rho(\gamma).\tilde s(x))=\delta_N(\tilde s(\gamma.x))=g.\delta_M(\gamma.x))=g.\delta_M(x)=\delta_N(\tilde s(x)) $$
We thus get that $\rho_N(\rho(\gamma))$ fixes pointwise an open subset of ${\bf S}^n$, hence is the identical transformation. We conclude that $\rho(\Gamma_M) \subset \text{Ker }\rho_N$, hence the map $\tilde s$ induces a conformal immersion $\sigma : M \to \mcn$, where $\mcn$ is the holonomy covering of $N$.  By construction, $\sigma$ is a lift of $s$.
\end{preuve}

\subsection{Cauchy completion of a conformally flat structure}
\label{sec.cauchy.completion}
The normal Cartan connection associated to a conformal structure allows to define an abstract notion of ``conformal boundary", derived from the $b$-boundary construction introduced in \cite{schmidt2}. We sketch the construction of this boundary below. More details are available in \cite[Sections 2 \text{and} 4]{frances.maximal}.  
 Fix once for all a basis $X_1, \ldots ,X_s$ of the Lie algebra $\lieg:=\lieo(1,n+1)$.  Given a Riemannian manifold $(M,g)$, with $\text{dim }M \geq 3$, let us call $(M,\hm,\om)$ the normal Cartan bundle associated to the conformal structure defined by $g$.  Denote by  ${\mathcal R}$ the frame field on $\hm$ defined by ${\mathcal R}(\hx)=((\omega_{\hx}^M)^{-1}(X_1), \ldots, (\omega_{\hx}^M)^{-1}(X_s))$.  This determines uniquely a Riemannian metric $\rho^M$ on $\hm$ having the property that ${\mathcal R}(\hx)$ is $\rho_{\hx}^M$-orthonormal for every $\hx \in \hm$.  The Riemannian metric $\rho^M$ defines a distance $d_M$ on $\hm$ by the formula 
 $$ d_M(\hx,\hy)=\frac{\delta_M(\hx,\hy)}{1+\delta_M(\hx,\hy)}$$
 where:
 
 - $\delta_M(\hx,\hy)$ is the infimum of the $\rho^M$-lengths of piecewise $C^1$ curves joining $\hx$ and $\hy$ if $\hx$ and $\hy$ are in the same connected component of $\hm$. 
 
 -   $\delta_M(\hx,\hy)=-2$ otherwise.
 
 One can look at the Cauchy completion $\hm_c$ of the metric space $(\hm, d_M)$, and define the Cauchy boundary $\partial_c\hm$ as $\partial_c\hm:=\hm_c \setminus \hm$.  Given $p \in P$, the right multiplication $R_p$ is Lipschitz with respect to $d_M$, and the right action of $P$ extends continuously to $\hm_c$.  The {\it conformal Cauchy completion} of $(M,g)$ is defined as the quotient space $M_c:=\hm_c\,/\,P$. 

Let us illustrate the construction in the case of the standard sphere ${\bf S}^n$, where the  conformal Cartan bundle is  identified with the Lie group $G=PO(1,n+1)$, and the Cartan connection is merely the Maurer-Cartan form $\omega^G$.  The Riemannian metric $\rho^G$ constructed as above is left-invariant on $G$, so that $(G,\rho^G)$ is a homogeneous Riemannian manifold, hence complete.  We infer that $G_c=\emptyset$, and the conformal Cauchy boundary of ${\bf S}^n$ is empty as well.

 Generally, the action of $P$ on $\hm_c$ is very bad behaved near points of $\partial_c\hm$, so that the space $M_c$ may not be Hausdorff.  It is thus quite remarkable that $M_c$ is Hausdorff when $(M,g)$ admits a conformal  immersion in the standard  sphere ${\bf S}^n$, as shows the following proposition.

\begin{proposition}
\label{prop.completion}
Let $M$ be a $n$-dimensional conformally flat manifold, $n\geq 3$. Assume there exists a conformal immersion $\delta: M \to {\bf S}^n$.   Then:
\begin{enumerate}
\item The conformal Cauchy completion $M_c$ is a Hausdorff space, in which $M$ is a dense open subset.
\item The conformal immersion $\delta$ extends to a continuous map $\delta: M_c \to {\bf S}^n$.
\item Every conformal diffeomorphism ${\varphi}$ of $M$ extends to a homeomorphism of $M_c$.
\end{enumerate}
\end{proposition}

\begin{preuve}
We call $\rho^M$ and $\rho^G$ the Riemannian metrics constructed on $\hm$ and $G$ as explained above, {\it using a same basis $X_1,\ldots,X_s$ of $\lieo(1,n+1)$}. The conformal immersion $\delta: M \to {\bf S}^n$ lifts to an isometric immersion $\hd : (\hm,{\rho}^M) \to (G,\rho^G)$. As a consequence, $\hd: (\hm,d_M) \to (G,d_G)$ is $1$-Lipschitz.  Because $(G,d_G)$ is a complete metric space, $\hd$ extends to a $1$-Lipschitz map $\hd: (\hm_c,d_M) \to (G,d_G)$.  This extended map $\hd$ is still $P$-equivariant for the (extended) action of $P$ on $\hm_c$ and on $G$. Every conformal diffeomorphism $\varphi \in \text{Conf}(M)$ lifts to an isometry $\hat \varphi$ of $(\hm,\rho^M)$, hence extends to an isometry, still denoted $\hat \varphi$ on $(\hm_c,d_M)$.  
The action of $P$ is free and proper on $\hm_c $ because the right action of $P$ on $G$ is free and proper, and $\hd$ maps $\hm_c $ continuously and $P$-equivariantly on $G$.  As a consequence, $M_c=\hm_c\,/\,P$ is Hausdorff. The map $\hd : \hm_c \to G$ induces a continuous $\delta: M_c \to G\,/\,P={\bf S}^n$, extending $\delta$.  Finally, for every $\varphi \in \text{Conf}(M)$, the homeomorphism $\hat \varphi: \hm_c \to \hm_c$ commutes with the right action of $P$, hence induces a homeomorphism $\varphi : M_c \to M_c$.  
\end{preuve}

\section{Proof of the local classification Theorem}
\label{sec.proof.local}
In this section, we  prove Theorem \ref{thm.local}. We are under the assumptions of the theorem, namely $s: L \setminus \Lambda \to N$ a conformal immersion, where $\Lambda$ is an essential singular set satisfying $\dhl=0$. We assume also that the singular set is essential and minimal in the sense that $\Lambda=\Lambda_{pole} \cup \Lambda_{ess}$, with $\Lambda_{ess}\not = \emptyset$. As explained in the introduction, because of Theorem \ref{thm.extension}, this hypothesis $\Lambda_{rem}=\emptyset$ is harmless.  By proposition \ref{prop.essential.flat} we know that both $L$ and $N$ are conformally flat manifolds.

\subsection{The target manifold $N$ is Kleinian}
\label{sec.target}
We call $\mcn$ the holonomy covering of $N$.  There is $\Gamma$ a discrete subgroup of conformal transformations of $\mcn$, acting freely properly discontinuously on $\mcn$ such that $N$ is conformally diffeomorphic to $\mcn\,/\,\Gamma$.  Showing that $N$ is Kleinian amounts  to show that $\mcn$ is conformally diffeomorphic to an open subset of ${\bf S}^n$.   The upshot of the  proof is as follows: we are going to construct a bigger $n$-dimensional conformal manifold $\mcn^{\prime}$, in which $\mcn$ embeds conformally as an open subset, and such that the action of $\Gamma$ extends conformally to $\mcn^{\prime}$.  The point is that the extended action of $\Gamma$ on $\mcn^{\prime}$  is no longer proper, what forces  $\mcn^{\prime}$ to be conformally equivalent to ${\bf S}^n$ or the Euclidean space (see Theorem \ref{thm.obata} below). Because $\mcn$ embeds conformally into $\mcn^{\prime}$, it is conformally diffeomorphic to an open subset of the sphere, as desired. 

\begin{theorem}[\cite{ferrand},\cite{schoen},\cite{frances1}]
\label{thm.obata}
Let $(M,g)$ be a Riemannian manifold of dimension  $n \geq 2$.  If the group of conformal transformations $\text{Conf}(M)$ does not act properly on $M$, then $(M,g)$ is conformally diffeomorphic to the standard sphere ${\bf S}^n$, or to the Euclidean space ${\bf R}^n$.
\end{theorem}
A version of the theorem for the identity component of the conformal group, and for compact manifolds, originally appeared in \cite{obata}.

We are now explaining how one can construct a manifold $\mcn^{\prime}$ with the properties listed  above.

In the remaining of this section, we  pick $x_{\infty} \in\Lambda_{ess}$, and $U$ a  connected neighborhood of $x_{\infty}$ in $L$, such that $U$ is conformally diffeomorphic to an open subset of the sphere ${\bf S}^n$.  Lemma \ref{lem.lift} ensures that the conformal immersion $s: U \setminus \Lambda \to N$ lifts to a conformal immersion $\si:U \setminus \Lambda \to \mcn$.  By definition of the holonomy covering, there exists a conformal immersion $\delta: \mcn \to {\bf S}^n$.  Then $\delta \circ \si: U \setminus \Lambda \to {\bf S}^n$ is a conformal immersion from $U \setminus \Lambda$ to an open subset of the sphere.  Because $\dhl=0$, $U \setminus \Lambda$ is a connected open subset of ${\bf S}^n$, and Liouville's theorem ensures that $\delta \circ \si$ is the restriction of a M\"obius transformation.  In particular, it is injective, and so is $\si$.  We thus get that $\si: U \setminus \Lambda \to \mcn$ is a conformal embedding.

In the following, we denote by $(\mcn,\hmcn,\omega^{\mcn})$ the normal Cartan bundle associated to the conformal structure on $\mcn$.  As in section \ref{sec.cauchy.completion}, we define the Riemannian metric $\rho^{\mcn}$ on $\hmcn$,  the associated distance $d_{\mcn}$,  $\hmcn_c$ the Cauchy completion of $(\hmcn,d_{\mcn})$, and $\mcn_c$ the conformal Cauchy completion of $\mcn$.   The distance on $\hmcn_c$ is still denoted $d_{\mcn}$.  

\begin{lemme}
\label{prop.point.infini}
The conformal embedding $\sigma: U \setminus \Lambda \to \mcn$ extends to a continuous map $\sigma: U \to \mcn_c$, which is a homeomorphism from $U$ onto an open subset $W \subset \mcn_c$.  The extendend map $\sigma$ sends $\Lambda \cap U$ into $\partial_c{\mcn}:=\mcn_c \setminus \mcn$.
\end{lemme}

\begin{preuve} Let us call $\hat U$ and $\hat \Lambda$  the inverse images  of $U$ and $\Lambda$ in  $\hat L$.  The conformal immersion $\si$ lifts to an isometric immersion $\hsi: (\hat U \setminus \hat \Lambda,\rho^{L}) \to (\hmcn,\rho^{\mcn})$.  Call $d_{U}$  (resp. $d_{U \setminus \Lambda}$) the distance induced by the Riemannian metric $\rho^L$ on the open set $\hat U$ (resp. $\hat U \setminus \hat \Lambda$).  Because $\hat\Lambda \cap \hat U$ has $(dim(\hat U)-1)$-dimensional Hausdorff measure  zero,  we get that $d_U=d_{U \setminus \Lambda}$ (this fact is probably standard. The reader can find a proof in \cite[Lemma 3.3]{frances.maximal}).  As a consequence, the map $\hsi: (\hat U \setminus \hat \Lambda,d_{U \setminus \Lambda}) \to (\hmcn,d_{\mcn})$, which is $1$-Lipschitz, is also $1$-Lipschitz if we put the metric $d_U$ on $\hat U \setminus \hat \Lambda$.   Hence, it  extends to a $1$-Lipschitz map $\hsi: (\hat U,d_U) \to (\hmcn_c,d_{\mcn})$.  This map is $P$-equivariant on the dense open subset $\hat U \setminus \hat \Lambda$, hence  on $\hat U$, and defines an extension of $\si$ to a continuous  map $\si: U \to \mcn_c$.

We are now going to show that the map $\si: U \to \mcn_c$ is open.  

Because  $\si: U \setminus \Lambda \to \mcn$ is an embedding, it is open on $U \setminus \Lambda$. It is thus enough to check that whenever $x \in\Lambda \cap U$, and $V \subset U$ is an open set containing $x$, the image $\sigma(V)$ is a neighborhood of  $z:=\sigma(x)$.  
Let $\hx \in \hat U$ be a point in the fiber of $x$, let $\hz=\hsi(\hx) \in \hmcn_c$, and let $r>0$ be very small so that $\overline{B(\hx,r)} $, the closure of the  ball of radius $r$ for $\rho^L$,  is compact and included in  $\hat V:=\pi_L^{-1}(V)$.   We claim that if $B(\hz,\frac{r}{5})$ denotes the metric ball centered at $\hz$ and of radius $\frac{r}{5}$ in $(\hmcn_c,d_{\mcn})$, we have the inclusion $B(\hz,\frac{r}{5}) \subset \hsi(\overline{B(\hx,r)})$, what will be enough to conclude, because the projections $\hat V \to V$ and $\hmcn_c \to \mcn_c$ are open maps.
Let us consider $\hz^{\prime} \in \hmcn_c$ such that $d_{\mcn}(\hz^{\prime},\hz)<\frac{r}{4}$.  Let us consider $(\hx_k)$ a sequence of $\hat U \setminus \hat \Lambda$ converging to $\hx$, and $(\hz_k^{\prime})$ a sequence of $\hmcn$ converging to $\hz^{\prime}$.  We consider indices $k$ large enough, so that {the points } $\hz_k:=\hsi(\hx_k)$  and $ \hz_k^{\prime}$ satisfy 
$$ d_{\mcn}(\hz_k,\hz_k^{\prime}) \leq \frac{r}{2}$$
and 
$$ d_U(\hx_k,\hx)<\frac{r}{5}.$$
  There is a curve $\beta_k: [0,1] \to \hmcn$ joining $\hz_k$ to $\hz_k^{\prime}$, and  having $\rho^{\mcn}$-length smaller than $\frac{3r}{4}$.  The key point is that there exist a lift $\alpha_k: [0,1] \to B(\hx,r) \setminus \hat \Lambda$, such that $\alpha_k(0)=\hx_k$ and $\hsi \circ \alpha_k=\beta_k$.  Let us see why it is true.  Let $t_{\infty}:=\sup \{ t \in [0,1], \text{ the lift } \alpha_k \text{ exists on } [0,t[   \}$.  Because $\hat \sigma: (\hat U \setminus \hat \Lambda,\rho^{L}) \to (\hmcn,\rho^{\mcn})$ is an isometric immersion, ${\alpha_k}_{|[0,t_{\infty}[}$ has finite length, so that $\hy_{\infty}:=lim_{t\to t_{\infty}} \alpha_k(t)$ exists.   Moreover, the $\rho^L$ length of ${\alpha_k}_{|[0,t_{\infty}[}$ is smaller than $\frac{3r}{4}$, so we get $d_U(\hx,y_{\infty})<r$, and $\hy_{\infty} \in B(\hx,r)$. If we prove that $\hy_{\infty} \not \in \hat \Lambda \cap B(\hx,r)$, we will get that $\alpha_k$ exists on $[0,1]$.  As we saw, the immersion $\si: U \setminus \Lambda \to \mcn$ is an embedding, so Theorem \ref{thm.extension.embeddings} ensures that all points of $\Lambda \cap U$ are either removable or poles with respect to $\sigma$.  Since $\sigma$ is a lift of $s$, any point of $\Lambda$ which is removable for $\sigma$ is removable  for $s$, and the minimality assumption on $\Lambda$ precisely says that there are no such points.  We conclude that every point of $\Lambda \cap U$ is a pole for $\sigma$. Hence, if we had  $\hy_{\infty}  \in \hat \Lambda \cap B(\hx,r)$, then $\hat \si (\alpha_k(t))$ should leave every compact subset of $\hmcn$ as $t \to t_{\infty}$, a contradiction with $ \beta_k([0,1]) \subset \hmcn$.
  
  The end point $\hx_k^{\prime}$ of $\alpha_k$ is mapped to $\hz_k^{\prime}$ by $\hsi$.  By compactness of $\overline{B(\hx,r)}$, we get a point $\hx^{\prime} \in \overline{B(\hx,r)}$ such that $\hsi(\hx^{\prime})=\hz^{\prime}$, what concludes the proof that $\si: U \to \mcn_c$ is open.
It remains to  check that it is injective to get that $\si$ maps $U$ homeomorphically onto its image $W$.  Let us assume for a contradiction that there are $x_1 \not =x_2$ in $U$ such that $\si(x_1)=\si(x_2)=y$.  Because $\si$ is open, there are $U_1$ and $U_2$ two disjoint open subsets of $U$ such that $\si(U_1) \cap \si(U_2)$ contains an open set $V$.  Now $\si(U_1 \setminus \Lambda) \cap V$ and $\si(U_2 \setminus \Lambda) \cap V$ being two dense open subsets of $V$, they must intersect, contradicting the injectivity of $\si$ on $U \setminus \Lambda$.  

We showed above that all points of $\Lambda \cap U$ are poles for the embedding $\si : U\setminus \Lambda \to \mcn$, what implies $\sigma(\Lambda) \subset \partial_c \mcn$. 
\end{preuve}

\begin{corollaire}
\label{coro.ouvert.sphere}
The holonomy covering $\mcn$ is conformally diffeomorphic to an open subset of ${\bf S}^n$, and $N$ is a Kleinian manifold.
\end{corollaire}

\begin{preuve}
Let us call $\rho : \Gamma \to PO(1,n+1)$ the group homomorphism satisfying the equivariance relation:
\begin{equation}
\label{eq.equi}
 \delta \circ \gamma=\rho(\gamma) \circ \delta
 \end{equation} 
  for every $\gamma \in \Gamma$.  We saw in Proposition \ref{prop.completion} that the action of $\Gamma$ extends to an action by homeomorphisms on $\mcn_c$, and that $\delta$ extends to a continuous map  $\delta: \mcn_c \to \Sn$.  In particular, by density of $\mcn$ in $\mcn_c$, the equivariance relation (\ref{eq.equi}) still holds on $\mcn_c$.  Let us define $\mcn^{\prime}:=\mcn \cup \bigcup_{\gamma \in \Gamma}\gamma(W)$. It is an open subset of $\mcn_c$, and in particular it is Hausdorff by Proposition \ref{prop.completion}.  By the previous proposition, the map $\delta \circ \si : U \to {\bf S}^n$ is continuous and coincides with the restriction of a M\"obius transformation on the dense open set $U \setminus \Lambda$.  Hence it is the restriction of a M\"obius transformation.  In particular $\delta: W \to {\bf S}^n$ is a homeomorphism on its image.  By  relation (\ref{eq.equi}), for every $\gamma \in \Gamma$, $\delta: \gamma(W) \to {\bf S}^n$ is a homeomorphism on its image as well.  From those remarks, we infer that $\mcn^{\prime}$ is a second countable Hausdorff space. The topological immersion $\delta: \mcn^{\prime} \to {\bf S}^n$ yields an atlas  which endows ${\mcn}^{\prime}$ with a structure of smooth conformally flat manifold, the conformal structure  ${\mathcal C}$ on $\mcn^{\prime}$ extending that of $\mcn$.  The equivariance relation (\ref{eq.equi}), available on ${\mcn}^{\prime}$, tells that in the charts of this atlas, the action of $\gamma \in \Gamma$ reads as the restriction of the action of $\rho(\gamma) \in PO(1,n+1)$.  In particular, $\Gamma$ acts as a group of smooth conformal transformations  of $(\mcn^{\prime},{\mathcal C})$.  
  
  We now use  Proposition \ref{prop.plat.injectif}, which implies that there exists a $G_{\delta}$-dense subset ${\mathcal G}$ of $N$ such that for every $y \in {\mathcal G}$, the fiber $s^{-1}\{y \}$ acumulates on our point $x_{\infty} \in \Lambda_{ess}$.  Because $\si$ is a lift of $s$, we get a sequence $(\gamma_k)$ of $\Gamma$, and a point $z \in \mcn$, such that $\gamma_k.z$ converges to $\si(x_{\infty})$.   We claim that the sequence $(\gamma_k)$ is not relatively compact in the conformal group of $\mcn^{\prime}$.   Indeed, if  it were not the case, $(\gamma_k)$ would preserve a Riemannian metric on $\mcn^{\prime}$, and then the function ``distance to $\partial_c\mcn \cap \mcn^{\prime}$" should be $\Gamma$-invariant.  Now, by  Lemma \ref{prop.point.infini}, $\sigma(x_{\infty}) \in \partial_c \mcn \cap \mcn^{\prime}$ so that the property $\gamma_k.z \to \si(x_{\infty})$ would lead to a contradiction.  
  
  Because $(\gamma_k)$ is not relatively compact,  the conformal group of $\mcn^{\prime}$ does not act properly, and  Theorem \ref{thm.obata}  ensures that $(\mcn^{\prime},{\mathcal C})$ is conformally equivalent to the standard $n$-sphere or the Euclidean $n$-space.  We infer that $\delta: \mcn \to \Sn$ is injective (Liouville's theorem), and $N$ is a Kleinian manifold.
\end{preuve}

\begin{remarque}
Actually, because the manifold $\mcn^{\prime}$ is conformally flat, we just need the conclusions of Theorem \ref{thm.obata} for conformally flat manifolds, and this result is actually much easier to prove than the general case.
\end{remarque}

\subsection{End of the proof of Theorem \ref{thm.local}}

 We keep the notations of Section \ref{sec.target}. Thanks to the work done there, we know that the developping map $\delta : \mcn \to \Sn$ is injective, so that $\delta$ is a conformal diffeomorphism between $\mcn$ and a connected open subset $\Omega \subset \Sn$.  Identifying $\Gamma$ with $\rho(\Gamma)$, we see $\Gamma$ as a Kleinian group in $PO(1,n+1)$, and get a commutative diagram 
 $$\xymatrix{
\mcn \ar@{>}[d] ^{\pi_{\mcn}} \ar@{>}[rr] ^{\delta}&  &  \Omega \ar@{>}[d] ^{\pi}\\
N \ar@{>}[rr] ^{\psi}&  & \Omega \,/\, \Gamma
}
$$
where $\psi$ is a conformal diffeomorphism.  We already noticed that $\Gamma$ does not act properly on ${\mcn}^{\prime}$, so that $\Gamma$ is infinite.

Let us pick $x_{\infty} \in \Lambda$, and a connected neighborhood $U$ of $x_{\infty}$ in $L$, which is conformally diffeomorphic so an open subset of the sphere.  By lemma \ref{lem.lift}, the conformal immersion $s : U \setminus \Lambda$ lifts to a conformal immersion $\sigma : U \setminus \Lambda \to {\mathcal N}$.  Liouville's theorem ensures that $\varphi:=\delta \circ \sigma$ extends to a conformal immersion  $\varphi : U \to \Sn$. Let us call $V:=\varphi(U)$. On $U \setminus \Lambda$, the  relation $\pi \circ \varphi=\psi \circ s$ holds, so that $\varphi$ yields a one-to-one correspondence between points of $\Lambda \cap U$ which are essential  (resp. poles) for $s$ to points of $\overline{\Omega} \cap V$ which are essential (resp. poles) for $\pi$.  By the discution of Section \ref{sec.examples.essential},  $\varphi$ maps $U \cap \Lambda$ to $V \cap \partial \Omega$, and  $U \cap \Lambda_{ess}$ to  $V \cap \Lambda(\Gamma)$. This completes the proof of Theorem \ref{thm.local}.

\subsection{Consequences of the local classification theorem}
\label{sec.consequences}
We are now proving Corollary \ref{coro.picard}, which derives some properties of the conformal singularities from theorem \ref{thm.local}. Our standing assumptions and notations are those of the corollary. 

We first explain why $\Lambda_{ess}$ is closed. Let us consider $(x_k)$ a sequence of $\Lambda_{ess}$ which converges to $x_{\infty} \in \Lambda$.  From Proposition \ref{prop.plat.injectif}, we know that $\text{Clust}(x_k)=N$ for all $k \in \NN$.  Hence, if we fix $y$ and $y^{\prime}$ two distinct points of $N$, one can  build two sequences $(y_k)$ and $(z_k)$ in $\lss$ which converge to $x_{\infty}$, such that $s(y_k) \to y$ and $s(z_k) \to y^{\prime}$. It follows that $x_{\infty} \in \Lambda_{ess}$.   Now, thanks to theorem \ref{thm.extension}, we extend $s$ to a conformal immersion $s^{\prime} : L \setminus (\Lambda_{ess} \cup \Lambda_{pole})$.   Theorem \ref{thm.local} implies that $N=\Omega/\Gamma$, for an infinite Kleinian group $\Gamma$.  It is a classical fact that the limit set $\Lambda(\Gamma)$ is either a perfect set, or has at most two points (\cite[Th.  2.3, p 43]{apanasov}). If we are in the former case, Theorem \ref{thm.local} ensures that $\Lambda_{ess}$ is perfect. If $\Lambda(\Gamma)$ has one or two points, then again by Theorem \ref{thm.local}, all the points of $\Lambda_{ess}$ are isolated. 

Assume that $\Lambda_{pole}$ is nonempty, and let us show that  the closure $\overline{\Lambda_{pole}}$ is $\Lambda_{pole} \cup \Lambda_{ess}$. By Theorem  \ref{thm.extension}, there is no harm assuming that $\Lambda=\Lambda_{pole} \cup \Lambda_{ess}$.   If $\Lambda_{ess}$ is empty, the claim is clear.   Assume now that  $\Lambda_{ess}$ is nonempty. It is enough to check that every point of $\Lambda_{ess}$ is in the closure of $\Lambda_{pole}$.  Recall that by Theorem \ref{thm.local}, for each $x_{\infty}$, there is a neighborhood $U$ of $x_{\infty}$ and a commutative diagram
 $$\xymatrix{
U \setminus \Lambda \ar@{>}[d] ^s \ar@{>}[rr] ^{\varphi}&  & V \setminus \partial \Omega \ar@{>}[d] ^{\pi}\\
N \ar@{>}[rr] ^{\psi}&  & \Omega \,/\, \Gamma
}
$$
Moreover, $\varphi(U \cap \Lambda)=V \cap \partial \Omega$ and  $\varphi(U \cap \Lambda_{ess})=V \cap \Lambda(\Gamma)$. We infer that $\partial \Omega \setminus \Lambda(\Gamma)$ is nonempty, and we are reduced to show that every point in $\Lambda(\Gamma)$ is accumulated by points in $\partial \Omega \setminus \Lambda(\Gamma)$. But this is clear, because  if $z \in \partial \Omega \setminus \Lambda(\Gamma)$, we will have $\Gamma.z \subset \partial \Omega \setminus \Lambda(\Gamma)$ and $\overline{\Gamma.z}=\Lambda(\Gamma) \cup \Gamma.z$.  

We assume now  that  $\Lambda$ is minimal essential. We want to show that if $x_{\infty} \in \Lambda_{ess}$ and if $U$ is any neighborhood of $x_{\infty}$ in $L$, then $s(U \setminus \Lambda)=N$.    By Theorem \ref{thm.local}, the manifold $N$ is Kleinian, conformally diffeomorphic to $\Omega/\Gamma$, for some infinite discrete  $\Gamma$.  For any $z \in \Omega$, the closure of $ \Gamma.z$ contains $\Lambda(\Gamma)$.  In particular, if $z_{\infty} \in \Lambda(\Gamma)$ and if $V$ is a neighborhood of $z_{\infty}$ in ${\bf S}^n$, then $\pi(V \setminus \Lambda(\Gamma))=\Omega/\Gamma$.  Theorem \ref{thm.local}  implies directly that $s(U \setminus \Lambda)=N$.


Finally, let us assume that $\Lambda$ is a discrete set containing at least one essential singular point.  Thanks to Theorem \ref{thm.extension}, we can assume that $\Lambda=\Lambda_{pole} \cup \Lambda_{ess}$.  The second point of the corollary implies that in the presence of essential singular points, $\Lambda_{pole}$  is not closed as soon as  it is nonempty.  Because $\Lambda$ is discrete, we infer that $\Lambda_{pole}$ must be empty.  If $\Gamma$ is the infinite Kleinian group such that $N=\Omega/\Gamma$, then $\Lambda(\Gamma)$ has one or two points (if not, $\Lambda_{ess}$ would be perfect), and Theorem \ref{thm.local} actually implies that $\Omega = {\bf S}^n \setminus \Lambda(\Gamma)$, otherwise $\Lambda$ would contain poles.  We infer that if $\Lambda(\Gamma)$ has one point, $\Gamma$ is a discrete subgroup of conformal transformations of $\R^n$ acting freely properly discontinuously on $\R^n$.  Then, one checks easily that $\Gamma$ is a discrete subgroup of Euclidean motions, and $N$ is a Euclidean manifold.  If $\Lambda(\Gamma)$ has two points, then $N$ is conformally diffeomorphic to a quotient of $\R^n \setminus \{ 0 \}$ by an infinite discrete group of conformal transformations, namely $N$ is a generalized Hopf manifold.

\section{Proof of Theorem \ref{thm.global}}
\label{sec.proof.global}

We are now considering  thin essential conformal singular sets on a compact manifold $L$.   This  compactness assumption on $L$ allows us to prove:

\begin{proposition}
\label{prop.revetement}
Let $(L,g)$ and $(N,h)$ be two connected $n$-dimensional Riemannian manifolds, $n \geq 3$.  Let $\Lambda \subset L$ be a closed subset such that $\dhl=0$, and $s:\lss \to N$ a conformal immersion. If $L$ is compact, and $\Lambda=\Lambda_{pole} \cup \Lambda_{ess}$, then  $s:\lss \to N$ is a covering map onto $N$.
\end{proposition}

\begin{preuve}
Let $\alpha : [0,1] \to N$ be a continuous path, let $x_0\in \lss$ such that $s(x_0)=\alpha(0)$. We want to show  the existence of $\gamma: [0,1] \to \lss$, a lift of $\alpha$ satisfying $\gamma(0)=x_0$.   If we can not lift $\alpha$, there exists $t_{\infty} \in [0,1[$, and  $\gamma: [0,t_{\infty}[ \to \lss$ a lift of $\alpha: [0,t_{\infty}[ \to N$, such that  $\gamma(0)=x_0$ and  $\gamma(t)$ leaves every compact subset of $\lss$ as $t$ tends to $t_{\infty}$.  By compactness of $L$, for every sequence $(t_k)$ tending to  $t_{\infty}$, the set $A$ of cluster values of $\gamma(t_k)$ in $L$ is nonempty and included in $\Lambda$.  Let $x_{\infty}$ be a point of $A$. Since $s(\gamma(t_k))$ tends to $\alpha(t_{\infty})$, we get $ x_{\infty} \not \in\Lambda_{pole}$.  Hence we should have $x_{\infty} \in \Lambda_{ess}$.  But this is not possible.  Indeed, if $x_{\infty} \in \Lambda_{ess}$, we first assume, considering a subsequence of $(t_k)$, that $\gamma(t_k)$ tends to $x_{\infty}$. Then we use Lemma \ref{lem.pas.equicontinu}, and get the existence of a sequence $(t_k^{\prime})$ in $[0,t_{\infty}[$, which converges to $t_{\infty}$, such that $\gamma(t_k^{\prime})$ converges to $x_{\infty}$, and such that $s(\gamma(t_k^{\prime}))$ converges to $y^{\prime} \in N$, with $y^{\prime} \not = \alpha(t_{\infty})$.  This contradicts the fact that $\gamma$ is a lift of $\alpha_{|[0,t_{\infty}[}$.
\end{preuve}

We are now under the hypotheses of Theorem \ref{thm.global}: the manifold $L$ is compact and the singular set is minimal essential, {\it i.e} $\Lambda=\Lambda_{ess} \cup \Lambda_{pole}$.  Moreover, we do the assumption ${\mathcal H}^{n-2}(\Lambda)=0$.  This assumption will be useful through the following result.

\begin{lemme}
\label{lem.simply}
Let $M$ be a connected, simply connected, $n$-dimensional Riemannian manifold, $n \geq 3$.  Assume that $E$ is a closed subset of $M$ satisfying ${\mathcal H}^{n-2}(E)=0$. Then $M \setminus E$ is still simply connected.
\end{lemme}

\begin{preuve}
In the proof we denote by $\overline{D}$  the closed unit disk in $\R^2$. 
We first observe that if $F$ is a closed subset of $\R^n$ satisfying ${\mathcal H}^{n-2}(F)=0$, and  if $f:\overline{D} \to \R^n$ is a $C^1$-map, then there exists a $C^1$-map $\tilde f : \overline{D} \to \R^n \setminus F$ agreeing with $f$ on $\partial D$.  Indeed, let $\mu$ be a smooth function on $\R^2$, which is positive on $D$ and vanishes identically  on the complement of $D$.  Let $D_m$ be the closed disk centered at the origin of radius $1-\frac{1}{m}$, $m \in \NN^*$.  The map $\varphi: D_m \times \R^n \to \R^n$ defined by $\varphi(x,u)=\frac{1}{\mu(x)}(f(x)-u)$ is locally Lipschitz, hence we must have ${\mathcal H}^{n}(\varphi(D_m \times F))=0$.  As a consequence, there exists $\xi \in \R^n$ which is in the complement of $\bigcap_{m \geq 1}\varphi(D_m \times F)$.  Then $\tilde f: \overline{D} \to \R^n \setminus E$ defined by  $\tilde f(x):=f(x)-\mu(x)\xi$ has the required properties.

Now, let us consider a loop $\gamma: [0,1] \to M \setminus E$.  Perturbing $\gamma$ into $M \setminus E$, we can assume that $\gamma$ is $C^1$. Since $M$ is simply connected, there exists a $C^1$-map $h: \overline{D} \to M$ such that $h(e^{2i\pi t})=\gamma(t)$, $\forall t \in [0,1]$.  Now, covering some $\overline{D}_m$ (for $m$ so big that $h(D \setminus \overline{D}_m) \subset M \setminus E$) by a finite number of small disks ${\mathcal D}_j$ such that $h({\mathcal D}_j)$ is contained in a chart of $M$, and applying to each $h: \overline{\mathcal D}_j \to M$ the observation made at the begining of the proof, we can construct $\tilde h: \overline{D} \to M \setminus E$ which agrees with $h$ on $\partial D$.  This yields that $\gamma$ is homotopically trivial  inside $M \setminus E$.
\end{preuve}

Theorem \ref{thm.local} ensures that  $(L,g)$ and $(N,h)$ are conformally flat, and that $N$ is actually conformally diffeomorphic, via a diffeomorphism  $\psi$, to a Kleinian manifold  $\Omega\,/\,\Gamma$.  From Theorem \ref{thm.local}, we also get that the  boundary $\partial \Omega$ satisfies ${\mathcal H}^{n-2}(\partial \Omega)=0$, hence $\Omega$ is simply connected by Lemma \ref{lem.simply}.    Let us call ${\tilde L}$ the conformal universal covering of $L$ and  denote by $\pi_{L}: {\tilde L} \to L$  the associated covering map. We call ${\tilde \Lambda}$ the inverse image of $\Lambda$ by $\pi_L$. Observe that ${\tilde L} \setminus {\tilde \Lambda}$ is simply connected by Lemma \ref{lem.simply}. By Proposition \ref{prop.revetement}, our conformal immersion $s:\lss \to N$ is a covering, hence it lifts to a conformal diffeomorphism $\sigma: {\tilde L} \setminus {\tilde \Lambda} \to \Omega$. In particular
\begin{equation}
\label{eq.equi0}
\pi \circ\sigma=\psi\circ s \circ \pi_L
\end{equation}

 Apply Theorem \ref{thm.extension.embeddings} to get that ${\sigma}^{-1}: \Omega \to {\tilde L}$ extends to a conformal diffeomorphism $\sigma^{-1}: \Omega^{\prime} \to {\tilde L}$, where $\Omega^{\prime} \subset \Sn$ is an open subset containing $\Omega$.  We denote again by $\sigma: {\tilde L} \to \Omega^{\prime}$ the inverse map. Observe that $\sigma(\tilde \Lambda)=\Omega^{\prime} \cap \partial \Omega$. The map $\sigma$ induces a homomorphism $\rho : \pi_1(L) \to PO(1,n+1)$ such that for every $\gamma \in \pi_1(L)$, the equivariance relation $\sigma \circ \gamma=\rho(\gamma) \circ \sigma$ holds. The group $\Gamma^{\prime}:=\rho(\pi_1(L))$ is a discrete subgroup of $PO(1,n+1)$ acting freely properly discontinuously on $\Omega^{\prime}$.     Let us call $\pi^{\prime}: \Omega^{\prime} \to \Omega^{\prime}/\Gamma^{\prime}$  the conformal covering map. There is  a conformal diffeomorphism $\varphi : L \to \Omega^{\prime}/ \Gamma^{\prime}$ such that
\begin{equation}
\label{eq.equi1}
\pi^{\prime} \circ \sigma =\varphi \circ \pi_L
\end{equation}
Let us check that $\Omega^{\prime}=\Omega(\Gamma^{\prime})$. If $\Gamma^{\prime}$ is finite, the compactness of $L$ leads to $\Omega^{\prime}=\Sn$.  If $\Gamma^{\prime}$ is infinite, one has $\Omega^{\prime} \subset \Omega(\Gamma^{\prime})$, since the action of $\Gamma^{\prime}$ is proper on $\Omega^{\prime}$. On the other hand, the compacity of $L$ forces the action of $\Gamma^{\prime}$ to be nonequicontinuous at each point of $\partial \Omega^{\prime}$, yielding the inclusion $\partial \Omega^{\prime} \subset \Lambda(\Gamma^{\prime})$.  In any case, we get that $\Omega^{\prime}=\Omega(\Gamma^{\prime})$, as claimed in Theorem \ref{thm.global}.

Now, observe that $\Gamma^{\prime} \subset \Gamma$, because  for every $\gamma \in \pi_1(L)$,  relation (\ref{eq.equi0}) leads to the identity $\pi \circ \rho(\gamma) =\pi$ on $\Omega$.  Hence, the identity map of $\Omega$ induces a covering map $s^{\prime} : \Omega/\Gamma^{\prime} \to \Omega/\Gamma$, satisfying for every $y \in \Omega$
\begin{equation}
\label{eq.equi2}
s^{\prime} \circ \pi^{\prime}(y)=\pi(y)
\end{equation}
Observe that if we define $\Lambda^{\prime}=\pi^{\prime}(\Omega^{\prime} \cap \partial \Omega)$, then $\Omega/\Gamma^{\prime}$ is merely $M(\Gamma^{\prime}) \setminus \Lambda^{\prime}$.
Relations (\ref{eq.equi0}),  (\ref{eq.equi1}) and  (\ref{eq.equi2}) lead to the commutative diagram
  
  $$\xymatrix{
L \setminus \Lambda \ar@{>}[d] ^s \ar@{>}[rr] ^{\varphi}&  & M(\Gamma^{\prime}) \setminus \Lambda^{\prime} \ar@{>}[d] ^{s^{\prime}}\\
N \ar@{>}[rr] ^{\psi}&  & \Omega \,/\, \Gamma
}
$$

  The diffeomorphism $\varphi$ maps $\Lambda$ to $\Lambda^{\prime}$ because $\sigma$ maps $\tilde{\Lambda}$ to $\partial \Omega \cap \Omega^{\prime}$.  Finally, it is easily checked that the essential singular points of $\Lambda^{\prime}$ for $s^{\prime}$ are the $\pi^{\prime}$-images of the essential singular points of $\partial \Omega \cap \Omega^{\prime}$ for $\pi$, namely the points of $\partial \Omega \cap \Omega^{\prime}$ which are in $\Lambda(\Gamma)$. This means $\Lambda(\Gamma) \cap \Omega^{\prime} \not = \emptyset$, hence  $\Lambda(\Gamma^{\prime}) \subsetneq \Lambda(\Gamma)$.

\section{Isolated essential singularities on compact manifolds}
\label{sec.isolated}

Our aim in this section is to understand completely the conformal singularities $s: \lss \to N$, where $N$ is a compact manifold and $\Lambda$ is a finite number of essential singular points.  It turns out that very few possibilities arise, and they are listed in Theorem \ref{thm.singularite.ponctuelle} below.  First of all, let us enumerate some examples.

\subsection{Euclidean singularities on the sphere}
Let us consider an {\it infinite}  dicrete subgroup $\Gamma \subset (\R_+^* \times O(n)) \ltimes \R^n$, acting freely properly discontinuously on $\R^n$.    One checks that for the action to be free, $\Gamma$ must actually be a subgroup of $O(n) \ltimes \R^n$.  The quotient manifold $N=\R^n\,/\,\Gamma$ is then a Euclidean manifold.  We  see $\Gamma$ as acting conformally on $\Sn \setminus \{\nu\}$, fixing $\nu$, and consider the covering map $s: \Sn \setminus \{\nu\} \to N$.  It is a conformal immersion, and because $\Gamma$ is infinite, we have  $\Lambda(\Gamma)=\{\nu\}$. Hence, as we already saw, $\nu$ is an essential singular point for $s$.  A conformal singularity  $s: \Sn \setminus \{  \nu \} \to N$ as described above  will be refered to as {\it Euclidean singularity on the sphere.}

\subsection{Singularities of Hopf type on the sphere}
\label{sec.hopf}
Let us now fix $o$ a second point on the sphere $\Sn$, distinct from the point $\nu$.  There is a conformal diffeomorphism mapping $\Sn \setminus \{o;\nu \}$ onto $\R^n \setminus \{0\}$. The group $G$  of conformal transformations of $\R^n \setminus \{0\}$ is generated by the inversion $\iota: x \mapsto -\frac{x}{||x||^2}$, and the group $\R_+^* \times O(n)$ of linear conformal transformations on $\R^n$.  Let us choose {\it an infinite} discrete group $\Gamma \subset G$  acting freely, properly and discontinuously on $\R^n \setminus \{0\}$.  It is not hard to check that $\Gamma$ has a finite index subgroup generated by a linear conformal contraction. As previously, the quotient $N=(\R^n \setminus \{0\})\,/\, \Gamma$ is called {\it a generalized Hopf manifold}.   The covering map $s:  \Sn \setminus \{o;\nu \} \to N$ is conformal, and because $\Gamma$ is infinite, both $\nu$ and $o$ are essential punctual singularities.  Conformal singularities $s: \Sn \setminus \{o; \nu \} \to N$ contructed as above   will be refered to as {\it  singularities of Hopf type on the sphere}.  

\subsection{Singularities of Hopf type on the projective space}


Let us go back to the previous construction, and assume that our infinite discrete subgroup $\Gamma \subset G$ contains the inversion $\iota$.  Then, the subgroup $\Gamma_o \subset \Gamma$ of transformations fixing individually the points $\nu$ and $o$ is normal in $\Gamma$.  Let us call $N_o$ the quotient manifold $(\R^n \setminus \{0\})\,/\, \Gamma_o$.  Because $\iota$ normalizes $\Gamma_o$, and because $\Gamma$ acts freely on $\R^n \setminus \{ 0\}$, $\iota$  induces a conformal involution $\overline{\iota}$ without fixed point on $N_o$.  The quotient $N_o/<\overline{\iota}>$ is actually conformally diffeomorphic to $N:=(\R^n \setminus \{0\})\,/\, \Gamma$.  The quotient of $ \Sn \setminus \{o;\nu \}$ by $<\iota>$ is conformally diffeomorphic to $\RP^n$ with a point $\nu$ removed. The natural covering map $\pi:  \Sn \setminus \{o;\nu \} \to N_o$ induces a conformal immersion $s: \RP^n\setminus \{ \nu \} \to   N$, for which $\nu$ is an essential singular point. Conformal singularities contructed in this way will be refered to as  {\it  singularities of Hopf type on  the projective space}.

\subsection{Classification result}

We are now investigating essential singular sets on compact manifolds, comprising only a finite number of points.  By Theorem \ref{thm.extension}, and the fourth point of Corollary \ref{coro.picard}, we just have to focus on the case where all the points are essential.  Then, it turns out  that the three kinds of singularities described in the previous section are the only possible.  

\begin{theorem} 
\label{thm.singularite.ponctuelle}
Let $(L,g)$ and $(N,h)$ be two connected $n$-dimensional Riemannian manifolds, $n\geq 3$, with $L$ compact.  Let $\Lambda:=\{p_1,\ldots,p_m\}$  be a finite number of points on $L$. Assume that  $s: L \setminus \Lambda \to N$ is a conformal immersion, such that each $p_i$ is an essential singular point for $s$.  Then $m=1$ or $m=2$ and:  

\begin{enumerate}
\item If $m=1$, either  there exists a Euclidean singularity on the sphere $s^{\prime}: \Sn \setminus \{\nu \} \to N^{\prime}$, a conformal diffeomorphism $\varphi : L \to \Sn$ sending $p_1$ to $\nu$ and a conformal diffeomorphism $\psi : N \to N^{\prime}$ making the diagram
$$\xymatrix{
L \setminus \{p_1\} \ar@{>}[d] ^s \ar@{>}[rr] ^{\varphi}&  & \Sn \setminus \{\nu \} \ar@{>}[d] ^{s^{\prime}}\\
N \ar@{>}[rr] ^{\psi}&  & N^{\prime}
}
$$
commute. 

 Or there exists a singularity of Hopf type on  the projective space $s^{\prime}: \RP^n \setminus \{\nu \} \to N^{\prime}$, a conformal diffeomorphism $\varphi : L \to \RP^n$ sending $p_1$ to $\nu $ and a conformal diffeomorphism $\psi : N \to N^{\prime}$ making the diagram
$$\xymatrix{
L \setminus \{p_1\} \ar@{>}[d] ^s \ar@{>}[rr] ^{\varphi}&  & \RP^n \setminus \{\nu \} \ar@{>}[d] ^{s^{\prime}}\\
N \ar@{>}[rr] ^{\psi}&  & N^{\prime}
}
$$
commute.
 
\item If $m=2$, there exists a singularity of Hopf type on the sphere  $s^{\prime}: \Sn \setminus \{o ; \nu \} \to N^{\prime}$, a conformal diffeomorphism $\varphi : L \to \Sn$ sending $\{ p_1;p_2\}$ to $\{ o; \nu \}$  and a conformal diffeomorphism $\psi : N \to N^{\prime}$ making the diagram
 $$\xymatrix{
L \setminus \{p_1;p_2\} \ar@{>}[d] ^s \ar@{>}[rr] ^{\varphi}&  & \Sn \setminus \{o; \nu \} \ar@{>}[d] ^{s^{\prime}}\\
N \ar@{>}[rr] ^{\psi}&  & N^{\prime}
}
$$

commute. 
  \end{enumerate}
\end{theorem}


\begin{preuve}
We first apply Theorem \ref{thm.local} in a neighborhood of any of the $p_i$'s.  We get that $N$ is conformally diffeomorphic to a Kleinian manifold $\Omega/\Gamma$, where the limit set $\Lambda(\Gamma)$ has one or two points (otherwise $\Lambda_{ess}$ would be a perfect set), and $\Omega=\Omega(\Gamma)$ (otherwise $\Lambda_{pole}$ would be nonempty).  

Assume first that $\Lambda(\Gamma)$ is made of a single point $\nu$.  The group $\Gamma$ is a discrete group of conformal transformations of $\Sn \setminus \{ \nu \}$, namely $\R^n$, which acts freely properly discontinuously on   $\R^n$.  As a consequence, $\Gamma$ is a discrete group of Euclidean motions, and $N$ is conformally diffeomorphic to  a Euclidean manifold $N^{\prime}=\R^n/\Gamma$. Theorem \ref{thm.global}  makes the structure of $L$ and $\Lambda$ explicit: there must be a subgroup $\Gamma^{\prime} \subset \Gamma$, with $\Lambda(\Gamma^{\prime}) \subsetneq \Lambda(\Gamma)$, as well as an open subset $\Omega^{\prime}$ properly containing $\Omega$, such that $L$ is conformally diffeomorphic to $\Omega^{\prime}/\Gamma^{\prime}$, and $\Lambda_{ess}$ is obtained as the quotient $(\Omega^{\prime} \cap \Lambda(\Gamma))/\Gamma^{\prime}$.  This implies in particular $\Lambda(\Gamma^{\prime})=\emptyset$, hence $\Gamma^{\prime}$ is finite, and because $\Gamma^{\prime}$ acts cocompactly on $\Omega^{\prime}$, we must have $\Omega^{\prime}=\Sn$.  Since the action of $\Gamma^{\prime}$ on $\Sn$ must be free, and $\Gamma^{\prime}$ fixes $\nu$, we infer  that $\Gamma^{\prime}$ is trivial.  We get   that $m=1$, $L$ is conformally diffeomorphic to $\Sn$, and we are in the first case of the theorem.

Assume now that $\Lambda(\Gamma)$ comprises two points $o$ and $\nu$.  Applying Theorem \ref{thm.global}, and with the same notations as above, we get that  $\Gamma$ is a discrete group in the conformal group of $\R^n \setminus \{  0\}$.  The limit set of the subgroup $\Gamma^{\prime}$ has two points or is empty, but because $\Lambda(\Gamma^{\prime}) \subsetneq \Lambda(\Gamma)$, we are in the second alternative: $\Gamma^{\prime}$ is once again finite, and $\Omega^{\prime}=\Sn$.  Because $\Gamma^{\prime}$ acts freely on $\Sn$, and leaves  $\{ o; \nu \}$ invariant, it is either trivial, or generated by a conformal involution of $\Sn$, without fixed point, and switching $o$ and $\nu$.  

It is not hard to check that  such a fixed-point free  involution switching $o$ and $\nu$  is conjugated, into the conformal group of $\R^n \setminus \{  0\}$, to the inversion $\iota: x \mapsto -\frac{x}{||x||^2}$, so if $\Gamma^{\prime}$ is nontrivial, there is no harm in assuming $\Gamma^{\prime}= <\iota>$.  Then $m=1$, $L$ conformally diffeomorphic to $\RP^n$, and we are in the second case of the theorem.

Finally, if $\Gamma^{\prime}$ is trivial, then $m=2$, $L$ is conformally diffeomorphic to $\Sn$ and we are in the third case of the theorem.
\end{preuve}

\section{Appendix: proof of Lemma \ref{dynamique-cones}}
\label{section.appendix}

We keep the notations introduced in section \ref{section.cartan}, especially sections \ref{section.exponential} and \ref{section.concormal.cones}.

We introduce $Q$ the quadratic form on $\R^{n+2}$ defined by $Q(x):=2x_0x_{n+1}+x_1^2+\ldots+x_n^2$.  We identify the group $G=O(1,n+1)$ with the group of linear transformations preserving $Q$, and we see $\Sn$ as the projectivization of the isotropic cone of $Q$.  We define $\nu:=[e_0]$ and $o:=[e_{n+1}]$ on $\Sn$.  The group $P$ is the stabilizer of $\nu$ in $PO(1,n+1)$.  Recall the splitting $\oo(1,n+1)=\lieg=\lien^- \oplus \R \oplus \oo(n) \oplus \lien^+$, where $\liep$ corresponds to $\R \oplus \oo(n) \oplus \lien^+$.  
The map $\rho: \lien^- \to {\bf S}^n \setminus \{ o\}$  defined  by $u \mapsto \exp_G(u).\nu$ is a diffeomorphism.
The map $j : {\bf R}^n \to {\bf S}^n$ defined by:
$$ (x_1, \ldots ,x_n) \mapsto [-\frac{Q(x)}{2},x_1, \ldots,x_n,1],$$
is a conformal diffeomorphism between the Euclidean space ${\bf R}^n$ and ${\bf S}^n \setminus \{\nu\}$.
In the sequel, we will call $\varphi$ the map $j^{-1} \circ \rho$. It is a diffeomorphism from $\lien^- \setminus \{ 0 \}$ ${\bf R}^n \setminus \{ 0 \}$ to ${\bf R}^n \setminus \{ 0 \}$.

We observe that $j$ intertwines the action of $P$ on ${\bf S}^n \setminus \{\nu\}$ and the affine action of $({\bf R}_+^* \times O(n)) \ltimes {\bf R}^n$ on ${\bf R}^n$. We will write the elements of $P$ in the affine form $\lambda A + T$, with $\lambda \in \R_+^*$, $A \in O(n)$, and $T \in \R^n$.

In the following, we will denote by $||.||$ the Euclidean norm on ${\bf R}^n$.  For a suitable choice of the $(\Ad O(n))-$ invariant scalar product $< \ , \ >$, $\varphi$ maps the Euclidean unit sphere on $S_{\lien^-}$.  It is then not hard to check that every conformal cone $C({\cal B},\lambda)$, with $\nu$ removed, is mapped by $j^{-1}$ to the set
$$\tilde{C}({\mathcal B},\lambda)\{ x=tu \in {\bf R}^n \ | \  t \in [\frac{1}{\lambda};\infty[\; , \ {u}\in {\varphi}({\cal B})  \}$$

Let $x \in {\bf R}^n$, and $u \in {\bf R}^n$  of Euclidean norm $1$. Then we define {\it the half-line} $[x,u)$ as the set
$$ [x,u):= \{ x+ tu \in {\bf R}^n \ | \ t \in {\bf R}_+ \}$$ 
 
The following lemma, the proof of which is left to the reader, will be useful in the sequel.
\begin{lemme}
\label{dynamique-segments}
Let $[x_k,v_k)$ be a sequence of half-lines in ${\bf R}^n$. Assume that whenever $v_{\infty}$ is a cluster value of $(v_k) $, then $-v_{\infty}$ is not a cluster value of   $\frac{x_k}{||x_k||}$.  Assume moreover that $x_k $ leaves every compact subset of $\R^n$. Then $[x_k,v_k)$ leaves every compact subset of ${\bf R}^n$.

\end{lemme}

We can now begin the proof of Lemma \ref{dynamique-cones}. Let us consider an unbounded sequence $(p_k)$ in $P$. Thanks to the chart $j$, we see $P$ as the conformal group  of ${\bf R}^n$.  The sequence $(p_k)$ then writes
 $$p_k : x \mapsto \lambda_k A_kx+\mu_k u_k,$$
  where $\lambda_k \in \Bbb R_+^*$,  $\mu_k \in \Bbb R_+$, $A_k \in O(n)$, and $||u_k|| =1$. Now, looking at a subsequence if necessary, we assume that $\lambda_k$, $\mu_k$, $\frac{\lambda_k}{\mu_k}$ all have limits in $\R_+^* \cup \{  +\infty\}$,  $u_k \to u_{\infty}$, and $A_k \to A_{\infty}$ in $O(n)$. The conclusions of Lemma \ref{dynamique-cones} won't be affected if we replace $p_k$ by $(A_k)^{-1}.p_k$, so that we may assume $p_k= \lambda_k Id + \mu_k u_k$.

Assume first that $\mu_k$ tends to $a \in \R_+$, and let $l_k$ be the translation of vector $-\mu_ku_k$. Clearly, $l_k \to l_{\infty}$ in $P$, whith $l_{\infty}$ the translation of vector $-au_{\infty}$, and $l_kp_k$ is just the homothetic transformation $x \mapsto \lambda_k x$, hence is in the factor ${\bf R}_+^* $ of $ P$. Since $(p_k)$ is unbounded, we can assume after taking a subsequence that $\lambda_k \to \infty$ or $\lambda_k \to 0$. In the first case, $l_kp_k.{\tilde C}({\cal B},\lambda)$ leaves every compact subset of ${\bf R}^n$, so that $l_kp_k.{ C}({\cal B},\lambda) \to \nu$, and we are in the first case  of the lemma. If $\lambda_k \to 0$, then $l_kp_k.{\tilde C}({\cal B},\lambda_k) \to {\tilde C}({\cal B},1)$, what yields  $(\Ad l_kp_k).{\cal C}({\cal B},\lambda_k) \to {\cal C}({\cal B},1)$, so that we are in the second case  of Lemma \ref{dynamique-cones}.  

We investigate now the case $\mu_k \to  \infty$, and $\frac{\lambda_k}{\mu_k} \to b_{\infty}$, $b_{\infty} \in \R_+$. Assume first that $b_{\infty}=0$, and let ${\cal B}^{\prime} \subset {\cal B}$ be a closed subball with nonzero radius, such that $-u_{\infty} \not \in {\varphi}({\cal B}^{\prime})$.   Let us consider a sequence of half-lines $[\frac{1}{\lambda}v_k,v_k)$ in ${\tilde C}({{\cal B}^{\prime}},\lambda)$. Here $(v_k)$ is a sequence of ${\varphi}({\mathcal B}^{\prime})$.  We observe that $p_k.[\frac{1}{\lambda}v_k,v_k)=[x_k,v_k)$, where $x_k=\frac{\lambda_k}{\lambda}v_k+\mu_ku_k$. 
Now, $\frac{x_k}{||x_k||}=\frac{\frac{\lambda_k}{\mu_k}\frac{1}{\lambda}v_k+u_k}{||\frac{\lambda_k}{\mu_k}\frac{1}{\lambda} v_k+u_k||}$, so that the only cluster value of $\frac{x_k}{||x_k||}$ is $u_{\infty}$.  We infer that if $v_{\infty}$ is a cluster value of $(v_k)$, then $-v_{\infty}$ can not be a cluster value of $\frac{x_k}{||x_k||}$. Writing $x_k=\mu_k(\frac{\lambda_k}{\lambda \mu_k}v_k+u_k)$, we check that $x_k \to \infty$.   Lemma \ref{dynamique-segments} ensures that $p_k.[\frac{1}{\lambda}v_k,v_k) \to \infty$. Since it is true for every sequence $[\frac{1}{\lambda}v_k,v_k)$, we get $p_k.{\tilde C}({{\cal B}^{\prime}},\lambda) \to \infty$. Hence $p_k.C({{\cal B}^{\prime}},\lambda) \to \nu$ and we are in the first case of Lemma \ref{dynamique-cones}.

If $b_{\infty} \not =0$, we choose ${\cal B}^{\prime} \subset {\cal B}$  a closed subball with nonzero radius, such that 
${\varphi}({\mathcal B}^{\prime}) \cap - {\varphi}({\mathcal B}^{\prime})= \emptyset$ and ${\varphi}({\mathcal B}^{\prime}) \cap \{ u_{\infty}; -u_{\infty} \} = \emptyset$. For $b$ near $b_{\infty}$, let us consider  the map
$$\psi: \alpha \mapsto \frac{x+\frac{u}{b \alpha}}{||x+\frac{u}{b \alpha}||}$$
One has $\psi(\alpha) \to \frac{x}{||x||}$ as $\alpha \to \infty$, and this uniformly with respect to $x \in {\varphi}({\mathcal B}^{\prime})$ and $(b,u)$ in a small compact neighborhood of $(b_{\infty},u_{\infty})$.  Hence, there is $\alpha_0 >0$, a $\delta>0$,  such that if $\alpha > \alpha_0$, $|b_{\infty}-b| \leq \delta$,  $||u-u_{\infty}|| \leq \delta$, and  $x \in {\varphi}({\mathcal B}^{\prime})$, then:
$$\sup_{v \in {\varphi}({\mathcal B}^{\prime})}|| \frac{b \alpha x + u}{||b \alpha x + u||} + v|| \geq \beta,$$
for some $\beta>0$.
Let us consider a sequence of half-lines $[2\alpha_0v_k,v_k)$ in ${\tilde C}({{\cal B}^{\prime}},\frac{1}{2 \alpha_0})$, where $v_k \in {\varphi}({\mathcal B}^{\prime})$.
We observe that $p_k.[2\alpha_0v_k,v_k)=[x_k,v_k)$, where $x_k=\lambda_k2\alpha_0v_k+\mu_ku_k$. 
Now, $\frac{x_k}{||x_k||}=\frac{\frac{\lambda_k}{\mu_k}2\alpha_0 v_k+u_k}{||\frac{\lambda_k}{\mu_k}2\alpha_0v_k+u_k||}$, so that for $k$ large
$$\sup_{v \in {\varphi}({\mathcal B}^{\prime})}|| \frac{x_k}{||x_k||} + v|| \geq \beta$$
It follows that if $v_{\infty}$ is a cluster value of $(v_k)$, then $-v_{\infty}$ can not be a cluster value of $\frac{x_k}{||x_k||}$.
Moreover, writing $x_k=\mu_k(\frac{\lambda_k}{\mu_k}2\alpha_0v_k+u_k)$, we see that because ${\varphi}({\mathcal B}^{\prime}) \cap \{u_{\infty}; -u_{\infty} \} = \emptyset$, $0$ is not a cluster value of $(\frac{\lambda_k}{\mu_k}2\alpha_0v_k+u_k)$, and $(x_k)$ tends to infinity.
We conclude thanks to lemma \ref{dynamique-segments} that $p_k.[2\alpha_0v_k,v_k) \to \infty$. Since it is true for every sequence $(v_k)$ of ${\varphi}({\mathcal B}^{\prime})$, we get $p_k.{\tilde C}({{\cal B}^{\prime}},\frac{1}{2\alpha_0}) \to \infty$, and we are in the first case of Lemma \ref{dynamique-cones}.

It remains to investigate the case where $\mu_k \to \infty$ and $\frac{\lambda_k}{\mu_k} \to \infty$. Let ${\cal B}^{\prime} \subset {\cal B}$ be a closed subball with nonzero radius, such that ${\varphi}({\mathcal B}^{\prime}) \cap -{\varphi}({\mathcal B}^{\prime}) = \emptyset$. Let $[\frac{1}{\lambda}v_k,v_k)$ be a sequence of half-lines in ${\tilde C}({{\cal B}^{\prime}},\lambda)$. For each integer $k$, $p_k.[\frac{1}{\lambda}v_k,v_k)=[x_k,v_k) $, with $x_k=\mu_k(\frac{\lambda_k}{\lambda \mu_k}v_k+u_k)$.  It is clear that $x_k \to \infty$.  The only cluster values of $\frac{x_k}{||x_k||}=\frac{\frac{1}{\lambda}v_k+\frac{\mu_k}{\lambda_k}u_k}{||\frac{1}{\lambda}v_k+\frac{\mu_k}{\lambda_k}u_k||}$ are included in ${\varphi}({\mathcal B}^{\prime})$.  We use once more lemma \ref{dynamique-segments} and conclude
$$p_k.{\tilde C}({{\cal B}^{\prime}},{\lambda}) \to \infty$$
We are again in the first case of Lemma \ref{dynamique-cones}.

Charles FRANCES\\
Laboratoire de Math\'ematiques, Bat. 425.\\
Universit\'e Paris-Sud 11.\\
91405 ORSAY.\\

\end{document}